\documentclass[12pt]{article}
\usepackage{unformatted}
\usepackage{url}
\date{February 7, 2017}
\usepackage{hyperref}
\hypersetup{
    colorlinks=true,
    linkcolor=blue,
    filecolor=magenta,
    urlcolor=cyan,
    citecolor=blue,
}
\usepackage[numbers]{natbib}
\usepackage{a4}
\usepackage{graphicx}
\usepackage{epstopdf}
\usepackage{ifthen}
\renewcommand{\cite}[2][nobrackets]{\ifthenelse{\equal{#1}{nobrackets}}{\Citet{#2}}{\Citep{#2}}}

\newcommand{\articletype}[1]{}
\newcommand{\received}[1]{}
\newcommand{\revised}[1]{}
\newcommand{\accepted}[1]{}
\newcommand{\corres}[1]{}

\usepackage{amssymb}
\usepackage{amsmath}
\usepackage{amsthm}

\usepackage{a4}
\usepackage{subeqnarray}

\usepackage{caption}
\usepackage{subcaption}
\usepackage{color}

\usepackage{pgfplots}
\pgfplotsset{compat=newest}
\usetikzlibrary{plotmarks}
\usepackage{grffile}

\usepackage{framed}
\usepackage[plain]{algorithm}
\usepackage{algpseudocode}

\journal{Journal of Computational Physics}

\newcommand{\be}{\begin{equation}}
\newcommand{\ee}{\end{equation}}
\newcommand{\bae}{\begin{eqnarray}}
\newcommand{\eae}{\end{eqnarray}}
\newcommand{\bse}{\begin{subeqnarray}}
\newcommand{\ese}{\end{subeqnarray}}
\newcommand{\brx}[1]{\left(#1\right)}

\definecolor{pink}{rgb}{1,0,1}
\definecolor{green}{rgb}{0,0.7,0}

\newcommand{\redii}[1]{\textcolor{black}{#1}}
\newcommand{\blueii}[1]{\textcolor{black}{#1}}
\newcommand{\pinkii}[1]{\textcolor{black}{#1}}

\newcommand{\red}[1]{\textcolor{black}{#1}}
\newcommand{\blue}[1]{\textcolor{black}{#1}}
\newcommand{\pink}[1]{\textcolor{black}{#1}}
\newcommand{\green}[1]{\textcolor{black}{#1}}

\newcommand{\Range}[1]{\redii{\textrm{Range}(#1)}}

\newcommand{\pd}[2]{\frac{\partial #1}{\partial #2}}
\newcommand{\pafg}[2]{\frac{\partial #1}{\partial #2}}

\newcommand{\secref}[1]{Section~\ref{sec:#1}}

\newcommand{\tabref}[1]{Table~\ref{tab:#1}}
\newcommand{\figref}[1]{Figure~\ref{fig:#1}}

\usepackage{xspace}
\newcommand{\methodname}{RAILS\xspace}
\newcommand{\methodnamelong}{Residual Approximation-based Iterative Lyapunov Solver\xspace}

\usepackage[pagewise]{lineno}
\newcommand{\llabel}[1]{}

\newcommand*\patchAmsMathEnvironmentForLineno[1]{%
  \expandafter\let\csname old#1\expandafter\endcsname\csname #1\endcsname
  \expandafter\let\csname oldend#1\expandafter\endcsname\csname end#1\endcsname
  \renewenvironment{#1}%
     {\linenomath\csname old#1\endcsname}%
     {\csname oldend#1\endcsname\endlinenomath}}%
\newcommand*\patchBothAmsMathEnvironmentsForLineno[1]{%
  \patchAmsMathEnvironmentForLineno{#1}%
  \patchAmsMathEnvironmentForLineno{#1*}}%
\AtBeginDocument{%
\patchBothAmsMathEnvironmentsForLineno{equation}%
\patchBothAmsMathEnvironmentsForLineno{align}%
\patchBothAmsMathEnvironmentsForLineno{flalign}%
\patchBothAmsMathEnvironmentsForLineno{alignat}%
\patchBothAmsMathEnvironmentsForLineno{gather}%
\patchBothAmsMathEnvironmentsForLineno{multline}%
\patchBothAmsMathEnvironmentsForLineno{eqnarray}%
\patchBothAmsMathEnvironmentsForLineno{subeqnarray}%
}

\def\txtd{{\textnormal{d}}}
\def\txte{{\textnormal{e}}}

\def\R{\mathbb{R}}
\def\E{\mathbb{E}}

\begin{document}

\begin{frontmatter}

\title{Continuation of Probability Density Functions using a Generalized 
Lyapunov Approach}

\author[rug]{S.~Baars}\ead{s.baars@rug.nl}
\author[cwi]{J.P.~Viebahn}\ead{viebahn@cwi.nl}
\author[uu]{T.E.~Mulder}\ead{t.e.mulder@uu.nl}
\author[tum]{C.~Kuehn}\ead{ckuehn@ma.tum.de}
\author[rug]{F.W.~Wubs}\ead{f.w.wubs@rug.nl}
\author[uu,cu]{H.A.~Dijkstra}\ead{h.a.dijkstra@uu.nl}

\address[rug]{Johann Bernoulli Institute for Mathematics and Computer Science, 
University of Groningen, P.O. Box 407, 9700 AK Groningen, The Netherlands}

\address[cwi]{Centrum Wiskunde \& Informatica (CWI), P.O. Box 94079, 1090 GB, Amsterdam, The Netherlands}

\address[uu]{Institute for Marine and Atmospheric research Utrecht, 
Department of Physics and Astronomy,  Utrecht University, 
Princetonplein 5, 3584 CC Utrecht, The Netherlands}

\address[tum]{Technical University of Munich,
Faculty of Mathematics, Boltzmannstr. 3, 85748 Garching bei M\"unchen, Germany}

\address[cu]{School of Chemical and Biomolecular Engineering, Cornell University, 
Ithaca, NY, USA}

\vfill
\begin{abstract}
  Techniques from numerical bifurcation theory are very useful to study transitions between steady fluid flow patterns and the instabilities involved.
  Here, we provide computational methodology to use parameter continuation in determining probability density functions of systems of stochastic partial differential equations near fixed points, under a small noise approximation.
  Key innovation is the efficient solution of a generalized Lyapunov equation using an iterative method involving low-rank approximations.
  We apply and illustrate the capabilities of the method using a problem in physical oceanography, i.e.~the occurrence of multiple steady states of the Atlantic Ocean circulation.
\end{abstract}


\end{frontmatter}

\newpage 
\section{Introduction}\label{sec:intro}
\pink{
\llabel{intro}Dynamical systems analysis of fluid flow phenomena has become a mature research direction \cite[]{Teman1988}.
The approach is rather complementary to direct numerical simulation of the governing equations in that the focus is \pinkii{on the direct computation of asymptotic sets in phase space (attracting forward or backward in time).}
The simplest of these sets are stable and unstable fixed points (steady states) and periodic orbits.
These sets play a major role in \pinkii{describing} transition behavior and associated pattern formation in many fluid flows, such as Rayleigh-B\'enard-Marangoni convection and the Taylor-Couette flow \cite[]{Koschmieder1993}.}

\pink{
Recently, an overview \cite[]{Dijkstra2014} was given of \pinkii{numerical techniques for computing} fixed points and periodic orbits for systems of partial differential equations (PDEs).
\pinkii{These PDEs are discretized} by any spectral, finite difference or finite element method.
Such discretizations result in systems of ordinary differential equations with algebraic constraints (e.g.~due to incompressibility) of high dimension, usually \pinkii{of} more than 10$^6$ degrees of freedom.
}
One of the most used methods is the pseudo-arclength continuation method \cite[]{Keller1977}, in which branches of steady states are parametrized by an arclength and the augmented problem (an equation representing the arc length normalization is added) is solved by the Newton-Raphson method.
The different methods are distinguished by whether the Jacobian matrix is explicitly computed (`matrix-based' methods) or whether matrix vector products of this matrix are used (`matrix-free' methods).
The applications shown in \cite{Dijkstra2014} \pinkii{range} from traditional flows \pinkii{and} free-surface flows to magneto-hydrodynamic flows and ocean flows.

A strong limitation of the numerical bifurcation techniques is that only relatively simple attracting sets can be computed.
Although tori can be determined in special cases \cite[]{Sanchez2013,Puigjaner2011}, the computation of fractal sets, often associated with chaotic dynamics of fluid systems, is out of the scope of these methods.
In this case one has to focus on ergodic properties of these systems, e.g.~on the determination of invariant measures, often represented by stationary probability density functions (PDFs).
Formally, one can determine these by solving the associated Liouville equation \pinkii{with a number of `spatial' dimensions equal to the degrees of freedom,} but solving these problems is not feasible for high-dimensional systems \cite[]{Lasota1994}.
Another limitation of numerical bifurcation methods is that when small-scale, non-resolved processes are represented as noise in the \pinkii{equations of the flow model}, one needs to solve the Fokker-Planck equation, which has even more computational obstacles than the Liouville equation \cite[]{Chang1970}.

In this paper, we focus on stochastic partial differential equations (SPDEs) describing fluid flows for which certain processes have been represented stochastically or for which the forcing of the flow has stochastic properties.
The direct way to investigate these flows is to use ensemble simulation techniques for many initial conditions \pinkii{to estimate the PDF for several observables of the flows} \cite[]{Slingo2011}.
Other methods which have been suggested use some form of model order reduction.
For example, a stochastic Galerkin technique forms the basis of the Dynamical Orthogonal Field method (DO) \cite[]{Sapsis2009}.
Non-Markovian reduced models can also be obtained by projecting on a basis of eigenvectors of the underlying deterministic system \cite[]{Chekroun2015B}.

In a deterministic-stochastic continuation method recently suggested by Kuehn \cite[]{Kuehn2012}, the results from fixed point computation in a deterministic model are used to obtain information on the stationary PDF of the stochastically forced system.
A restriction is that linearized dynamics near the fixed point adequately \pinkii{describe} the behavior of the stochastic system.
In this case, one only needs the leading-order linear approximation and determines the covariance matrix from the solution of a Lyapunov equation.
This provides all the information to determine the probability of sample paths.
The nice aspect of this method is that one can combine it easily with deterministic pseudo-arclength continuation methods.
However, in order to apply the approach to systems of PDEs with algebraic constraints, generalized Lyapunov equations have to be solved.

Direct methods to solve a generalized Lyapunov equation such as Bartels-Stewart algorithm \cite[]{Bartels1972} are based on dense matrix solvers and hence inapplicable for large systems.
Other existing methods which use low-rank approximations such as Extended and Rational Krylov subspace methods \cite[]{Simoncini2007,Druskin2011,Stykel2012,Druskin2014} and alternating directions implicit (ADI) based iterative methods \cite[]{Kleinman1968, Penzl1999} might also become expensive for large-dimensional problems, particularly \pinkii{when} trying to use previous initial guesses along a continuation branch.

\green{
The aim of this paper is to present new methodology to efficiently trace PDFs of SPDEs with algebraic constraints in parameter space.
In \secref{methods}, an extension of the approach suggested in \cite{Kuehn2012} and the novel procedure to efficiently solve a generalized Lyapunov equation numerically are presented.
In fact, this solves an open conjecture in \cite{Kuehn2015}, which states that it should be possible to reuse previous solutions of a continuation in a specialized Lyapunov solver.
We describe our application in \secref{test}, which is a model of the Atlantic Ocean circulation.
The numerical aspects and capabilities of the novel method for this application are shown in \secref{results}.
In \secref{conc}, we provide a summary and discuss the results.
}
\section{Methods}\label{sec:methods}

Any ocean-climate model consists of a set of conservation laws (momentum, mass, heat and salt), which are formulated as a set of coupled partial differential equations, that can be written in general form as \cite[]{Griffies2004}
\be
{\cal M}({\bf p}) \pafg{{\bf u}}{t}  = {\cal L}({\bf p}) {\bf u} +  {\cal N}({\bf u}, {\bf p}) 
+  {\cal F}({\bf u}, {\bf p}), 
\label{e:gen_sys}
\ee
where ${\cal L}$, ${\cal M}$ are linear operators, ${\cal N}$ is a nonlinear operator, $\bf u$ is the state vector, $\cal F$ contains the forcing of the system and ${\bf p}\in\R^{n_p}$ indicates a vector of parameters.
Appropriate boundary and initial conditions have to be added to this set of equations for a well-posed problem.

\subsection{Formulation of the problem}
When Eq.~(\ref{e:gen_sys}) is discretized, eventually a set of ordinary differential
equations with algebraic constraints \pinkii{arises,} which can be written as 
\be
M ({\bf p}) \frac{\txtd {\bf x}}{\txtd t}  = L({\bf p}) {\bf x} + N({\bf x}, {\bf p}) +  F({\bf x}, {\bf p}),  
\label{e:gen_sys_discrete}
\ee
where ${\bf x}\in\mathbb{R}^n$ is the state vector, \llabel{mass}$M({\bf p})\in\R^{n\times n}$ is a singular matrix of which every zero row is associated with an algebraic constraint, $L({\bf p})\in\R^{n\times n}$ is the discretized version of $\mathcal{L}$, and $F:\R^{n_p}\rightarrow\R^n$ and $N:\R^n\times \R^{n_p}\rightarrow\R^n$ are the finite-dimensional versions of the forcing and the nonlinearity respectively.
When noise is added to the forcing, the evolution of the flow can generally be described by a stochastic differential-algebraic equation (SDAE) of the form
\bae
\label{SDAE}
M({\bf p}) ~\txtd{\bf x}_t={\bf f}({\bf x}_t;{\bf p})~\txtd t+{\bf g}({\bf x}_t;{\bf p})~\txtd{\bf W}_t,
\eae
where ${\bf f}({\bf x}_t;{\bf p}) = L({\bf p}) {\bf x}_{t} + N({\bf x}_{t}, {\bf p}) +  F({\bf x}_{t}, {\bf p})$ is 
the right-hand side of Eq.~\eqref{e:gen_sys_discrete}, ${\bf W}_t\in \R^{n_w}$ is a vector of $n_w$-independent
standard Brownian motions~\citep{Gardiner2009}, and ${\bf g}({\bf x}_t;{\bf p})\in \R^{n_w\times n}$.

Suppose that the deterministic part of Eq.~\eqref{SDAE} has a stable fixed point ${\bf x}^*={\bf x}^*({\bf p})$ for a given range of parameter values.
Then linearization around the deterministic steady state yields \citep{Kuehn2012}
\bae\label{SDAElin}
M({\bf p}) ~\txtd{\bf X}_t=A({\bf x}^*;{\bf p}){\bf X}_t~\txtd t+B({\bf x}^*;{\bf p})~\txtd{\bf W}_t,
\eae
where $A({\bf x};{\bf p})\equiv(D_{\bf x}{\bf f})({\bf x};{\bf p})$ is the Jacobian matrix and 
$B({\bf x};{\bf p}) = {\bf g}({\bf x}^*;{\bf p})$.
From now on, we drop the arguments of the matrices $A, B$ and $M$.

In the special case that $M$ is a non-singular matrix, the equation~(\ref{SDAElin}) can be rewritten as
\bae\label{SDAElin2}
\txtd {\bf X}_t=M^{-1}A {\bf X}_t~\txtd t+M^{-1} B ~\txtd{\bf W}_t,
\eae
\blue {
which represents an $n$-dimensional Ornstein-Uhlenbeck (OU) process. \llabel{uo}
}
The corresponding stationary covariance matrix $C$ is determined from the following \emph{Lyapunov equation} \citep{Gardiner2009}
\bae\label{lyapunov}
M^{-1}AC+CA^{\top}M^{-\top}+M^{-1}BB^{\top}M^{-\top}=0. 
\eae
This equation can be rewritten as a \emph{generalized Lyapunov equation}
\bae\label{glyapunov}
ACM^{\top}+MCA^{\top}+BB^{\top}=0.
\eae
If $M$ is a singular diagonal matrix then Eq.~(\ref{SDAElin2}) does not apply.
However, if the stochastic part is non-zero only on the part where $M$ is non-singular (which occurs often when the noise is in the forcing of the flow), then Eq.~\eqref{SDAElin} can be written as
\bae\label{SDAEsystem}
\begin{pmatrix}0&0\\0&M_{22}\end{pmatrix}
\begin{pmatrix}\txtd{\bf X}_{t,1}\\ \txtd {\bf X}_{t,2}\end{pmatrix}
=
\begin{pmatrix}A_{11} & A_{12} \\ A_{21} & A_{22}\end{pmatrix}
\begin{pmatrix}{\bf X}_{t,1}\\{\bf X}_{t,2}\end{pmatrix}
\txtd t
+
\begin{pmatrix}0\\B_{2}\end{pmatrix}\txtd{\bf W}_t,
\eae
where $M_{22}$ represents the non-singular part of $M$.
Consequently, for non-singular $A_{11}$ we can separate Eq.~\eqref{SDAEsystem} into an algebraic part and an explicitly time-dependent part,
\bse\label{SDAEsystem2}
A_{11}{\bf X}_{t,1}+A_{12}{\bf X}_{t,2} &=& 0 ,\slabel{SDAEsystem2a}\\
M_{22}~\txtd{\bf X}_{t,2} &=& S{\bf X}_{t,2}~\txtd t+B_{2}~\txtd {\bf W}_t,\slabel{SDAEsystem2b}
\ese
where $S=A_{22}-A_{21}A_{11}^{-1}A_{12}$ is the Schur complement of $A$, which is well-defined since $A_{11}$ is nonsingular as the Jacobian $A$ has eigenvalues strictly in the left-half complex plane by the \llabel{ass}ass\blue{umption on deterministic stability of ${\bf x^*}$}.
Since $M_{22}$ is non-singular by construction we can find the stationary covariance matrix $C_{22}$ by solving the corresponding generalized Lyapunov equation
\bae\label{schurlyapunov}
SC_{22}M_{22}^{\top}+M_{22}C_{22}S^{\top}+B_{2}B_{2}^{\top}=0.
\eae

We remark that Eq.~\eqref{schurlyapunov} could alternatively be derived using an epsilon-embedding approach for the differential algebraic equations.
In order to find the full covariance matrix we use $C_{ij}=\E[{\bf X}_{t,i}{\bf X}_{t,j}^{\top}]$ together with Eq.~\eqref{SDAEsystem2a} which gives
\begin{align*}
  C_{12}&=-A_{11}^{-1}A_{12}\E[{\bf X}_{t,2}{\bf X}_{t,2}^{\top}]\\
  &=-A_{11}^{-1}A_{12}C_{22},\\
  C_{21}&=-\E[{\bf X}_{t,2}{\bf X}_{t,2}^{\top}]A_{12}^{\top}A_{11}^{-\top}\\
  &=-C_{22}A_{12}^{\top}A_{11}^{-\top}=C_{12}^\top,\\
  C_{11}&=A_{11}^{-1}A_{12}\E[{\bf X}_{t,2}{\bf X}_{t,2}^{\top}]A_{12}^{\top}A_{11}^{-\top}\\
  &=-A_{11}^{-1}A_{12}C_{21}.
\end{align*}

In summary, we can obtain an estimate of the covariance matrix $C$ of the stochastic dynamical system~\eqref{SDAE} by providing the matrices $A$, $B$ and $M$ and solving the corresponding generalized Lyapunov equation~\eqref{glyapunov}, or \eqref{schurlyapunov} in case $M$ is singular.
Once the covariance matrix $C$ is computed, the stationary PDF of the approximating OU-process, indicated by $p({\bf x})$ follows as \cite[]{Gardiner2009, Kuehn2011}
\bae
p({\bf x}; {\bf x}^*) =  \frac{1}{(2 \pi)^{\frac{n}{2}}} \mid C \mid^{-1/2} 
\txte^{-\frac{1}{2}({\bf x}-{\bf x}^*)^{\top} C^{-1} ({\bf x}-{\bf x}^*)} .
\eae
The limitations of this approach are, firstly, that only a PDF estimate is provided valid on 
a subexponential time scale before large deviations occur and, secondly, that only the local 
behavior near the steady state and Gaussian stochastic behavior of the system are obtained.

\subsection{A novel iterative generalized Lyapunov solver}
The type of systems of the form \eqref{glyapunov} that we want to solve are typically sparse and have a dimension $n = {\cal O}(10^5)$ or larger.
Solving systems of this size results in a $C$ that is generally a dense matrix of the same size.
This is computationally very expensive in terms of both time and memory.
Consequently, one cannot aim to compute the full $C$ but only a low-rank approximation of the form $C\approx VTV^{\top}$.
\blue{
In existing iterative solution methods for low-rank approximations \cite{Kleinman1968,Penzl1999,Saad1990,Simoncini2007,Stykel2012} the matrix $V$ is usually computed using repetitive products with $B$ in every iteration, for instance in such a way that it spans the Krylov spaces $\mathcal{K}_m(A,B)$ or $\mathcal{K}_m(A^{-1},B)$.
\llabel{bvec}In practice, however, the matrix $B$ might have many columns, which means that in every iteration of such a method many matrix-vector products have to be performed or many linear systems have to be solved.
These operations take up by far the largest amount of time in every iteration, which is why we would like a method that does not expand the search space with the same amount of vectors as the amount of columns in $B$.
}

\llabel{saad}The solution method we propose is based on a Galerkin projection, and is very similar to the method in \cite{Saad1990}.
It works by solving projected systems of the form
\begin{align*}
  V^{\top}AVTV^{\top}M^{\top}V+V^{\top}MVTV^{\top}A^{\top}V
	+V^{\top}BB^{\top}V&=0,
\end{align*}
where $C$ is approximated by a low-rank approximation $\tilde C=VTV^{\top}$.
Now if we take $\tilde A = V^{\top}AV$, $\tilde M=V^{\top}MV$, $\tilde B = V^{\top}B$, we get the smaller generalized Lyapunov equation
\begin{align}\label{proj_sys}
  \tilde A T \tilde M^{\top}+\tilde M T \tilde A^{\top} + \tilde B \tilde B^{\top} = 0, 
\end{align}
which can be solved by a dense solution routine \cite[]{Bartels1972}. 

A problem that arises when solving generalized Lyapunov equations in an iterative manner 
is computing an estimate of the residual
\begin{align}\label{resid}
  R=A\tilde CM^{\top}+M\tilde CA^{\top}+BB^{\top} , 
\end{align}
for some approximate solution $\tilde C$.
The (matrix) norm of this residual, which is generally a dense matrix, can be used in a stopping criterion.
The 2-norm of the residual matrix is equal to its spectral radius which is defined by the absolute value of the largest eigenvalue.
Since the residual is symmetric, approximations of the largest eigenpairs  can be computed using only a few steps of the Lanczos method \cite{Lanczos1950}.
Even though we cannot compute $R$ explicitly, it is possible to apply the Lanczos method to determine the eigenpair because only matrix vector products $R x$ are needed, which are evaluated as
\begin{align*}
\redii{Rx = A ( V ( T ( V^{\top} (M^{\top} x)))) + M (V ( T ( V^{\top} ( A^{\top} x)))) + B ( B^{\top} x) . }
\end{align*}
\blue {
The goal of our method is to compute the matrix $V$, which we could also view as a search space by considering its columns as basis vectors for a linear subspace of $\R^n$.
\llabel{orth}From now on we assume $V$ to be orthonormalized.
We suggest to expand $V$ in every iteration by the \redii{eigenvectors associated with the largest eigenvalues} of the residual, which we already obtained when computing the norm of the residual.
The reasoning behind this will be explained below.
}
\red {
\llabel{res}Now in every iteration, we solve the projected system \eqref{proj_sys}, but because we expanded our search space (with the largest components of the residual), we hope that the new residual is smaller.
}
The resulting algorithm for solving generalized Lyapunov equations is shown in Algorithm~\ref{alg:proj}.
Because of the peculiar choice of vectors to expand our space, we call this method the \methodnamelong (\methodname).
\begin{algorithm}[H]
\centering
\parbox{.95\textwidth}{
\begin{framed}
\begin{tabular}{lll}
        \textbf{input:}  & $A,B,M$ & The problem where $ACM^{\top}+MCA^{\top}+BB^{\top}=0$.\\
        & $V_1$ & Initial space.\\
        & $m$ & Dimension increase of the space  per iteration.\\
        & $l$ & Maximum amount of iterations.\\
        & $\epsilon$ & Convergence tolerance.\\
        \textbf{output:} & $V_k, T_k$& Approximate solution, where $C\approx V_kT_kV_k^{\top}$.
\end{tabular}
\begin{algorithmic}[1]
  \vspace{1mm}
  \State Orthonormalize $V_1$
  \State Compute $\tilde A_1=V_1^{\top}AV_1$
  \State Compute $\tilde M_1=V_1^{\top}MV_1$
  \State Compute $\tilde B_1=V_1^{\top}B$
  \For{$j= 1,\ldots,l$}
  \State Obtain $\tilde A_j=V_j^{\top}AV_j$ by only computing new parts
  \State Obtain $\tilde M_j=V_j^{\top}MV_j$ by only computing new parts
  \State Obtain $\tilde B_j=V_j^{\top}B$ by only computing new parts
  \State Solve $\tilde A_j T_j \tilde M_j^{\top} + \tilde M_j T_j \tilde A_j^{\top} + \tilde B_j \tilde B_j^{\top}=0$
  \State Compute the approximate largest $m$ eigenpairs ($\lambda_p$, $r_p$) of the residual $R_j$ using Lanczos
  \State Stop if the approximated largest eigenvalue is smaller than $\epsilon$ \label{alg:stop}
  \State $V_{j+1} = [V_j, r_1,\ldots,r_m]$
  \State Re-orthonormalize $V_{j+1}$
  \EndFor
\end{algorithmic}
\end{framed}}
\caption{\methodname algorithm for the projection based method for solving generalized Lyapunov equations.}
\label{alg:proj}
\end{algorithm}

\newtheorem{theorem}{Theorem}[section]
\newtheorem{proposition}{Proposition}[section]
\newtheorem{remark}{Remark}[section]
\newtheorem{corollary}{Corollary}[section]
\newcommand{\orth}[1]{\textrm{orth}(#1)}

\subsection{\green{Convergence analysis}}\llabel{props}
We will now show why we choose the \redii{eigenvectors associated with the largest eigenvalues} of the residual. Here we use $\orth{B}$ to denote the orthonormalization of $B$. We first show that in a special case, $V_k$ spans the Krylov subspace
\begin{align*}
  \mathcal{K}_k(A,B) = \{B, AB, \ldots, A^{k-1}B\}.
\end{align*}
\begin{proposition}\label{t1}
  If $M=I$, $B \in \mathbb{R}^{n\times m}$, where $m$ is also the amount of vectors we use to expand the space $V_k$ in every iteration and $V_1 = \orth{B}$, then $\Range{V_k} \subseteq \mathcal{K}_k(A,B)$.
  \begin{proof}
    For $k=1$ this is true by \redii{assumption}.
    Now say that in step $k$, $(\lambda,q)$ is an eigenpair of the residual $R_k$, and assume that $\Range{V_k} \subseteq \mathcal{K}_k(A,B)$.
    Then we can write
    \begin{align*}
      R_kq=\lambda q=AV_kq_1+V_kq_2+Bq_3
    \end{align*}
    where $q_1 = T_kV_k^{\top}q$, $q_2 = T_kV_k^{\top}A^{\top}q$, $q_3 = B^{\top}q$.
    From this it is easy to see that if we orthonormalize $q$ with respect to $V_k$, it is only nonzero in the direction of $AV_k$.
    Now we take $V_{k+1} = [V_k, Q_k]$, where \redii{the columns of} $Q_k$ are the eigenvectors \redii{associated with} the $m$ largest eigenvalues of $R_k$ orthonormalized with respect to $V_k$. Then
    \begin{align*}
      \Range{V_{k+1}} &\subseteq \Range{V_k} \cup \Range{AV_k}\\
      &\subseteq \mathcal{K}_k(A,B) \cup A\mathcal{K}_k(A,B) = \mathcal{K}_{k+1}(A,B).
    \end{align*}
  \end{proof}
\end{proposition}
\begin{remark}\label{r1}
\redii{
In case $Q_i$ has full rank for every $i=1,\ldots,k$ it is clear that actually the equality $\Range{V_k} = \mathcal{K}_k(A,B)$ holds in Proposition~\ref{t1}. This is what we observe in practice.
}
\end{remark}
\redii {
\llabel{saad2}From Proposition~\ref{t1} and Remark~\ref{r1}, we see that when we choose $m$ equal to the amount of columns of $B$, \methodname is equivalent to the method in \cite{Saad1990} as long as $Q_k$ has full rank in every iteration.
What is important, is that we want to take $m$ much smaller, in which case we assume that \redii{eigenvectors associated with the largest eigenvalues} of the residual $R_k$ point into the direction of the most important components of $AV_k$.
A similar result holds when we want to look in the Krylov space $\mathcal{K}_k(A^{-1},A^{-1}B)$.
}
\begin{proposition}\label{t2}
  If $M=I$, $B \in \mathbb{R}^{n\times m}$, where $m$ is also the amount of vectors we use to expand the space $V_k$ in every iteration, $V_1 = \orth{A^{-1}B}$ and $V_{k+1}=[V_k, Q_k]$, where $Q_k$ are the $m$ \redii{eigenvectors associated with the largest eigenvalues} of $A^{-1}R_k$ orthonormalized with respect to $V_k$, then $\Range{V_k} \subseteq \mathcal{K}_k(A^{-1},A^{-1}B)$.
  \begin{proof}
    This can be proved analogously to Proposition \ref{t1}.
  \end{proof}
\end{proposition}
\pink {
We show this result since most other iterative Lyapunov solvers include an operation with $A^{-1}$. \pinkii{An example} is the Extended Krylov method, which looks in the Krylov space \redii{$\mathcal{K}_{2k}(A,A^{-k}B)$.}
\pinkii{We remark} that our method, when we start with $V_1=\orth{[B, A^{-1}B]}$ and expand with $[Q_k, A^{-1}Q_k]$ is not equivalent to the Extended Krylov method.
}

\pink {
We know that if $V_k$ has $n$ orthogonal columns, it spans the whole space, so the solution is in there.
\pinkii{To show that our method has finite termination, we argue that the method has converged when the vectors we generate do not have a component perpendicular to $V_k$, which means that the size of the search space does not increase anymore.}
}
\begin{proposition}\label{t3}
  Take $M=I$ and $B \in \mathbb{R}^{n\times m}$. After $k$ steps of \methodname, the residual $R_k$ has an eigenpair \redii{$(\lambda,q)$ with $\lambda = \Vert R_k\Vert_2$}. If $q \in \Range{V_k}$, then $V_kT_kV_k^{\top}$ is the exact solution.
  \begin{proof}
    We have
    \begin{align*}
      R_kq=\lambda q=AV_kT_kV_k^{\top}q+V_kT_kV_k^{\top}A^{\top}q+BB^{\top}q.
    \end{align*}
    Since $q \in \Range{V_k}$ and $V_k$ is orthonormalized, it holds that $q=V_kV_k^{\top}q$. So then
    \begin{align*}
      \lambda q=\lambda V_kV_k^{\top}q &= V_kV_k^{\top}AV_kT_kV_k^{\top}q+V_kT_kV_k^{\top}A^{\top}V_kV_k^{\top}q+V_kV_k^{\top}BB^{\top}V_kV_k^{\top}q\\
                                       &= V_k\tilde A_kT_kV_k^{\top}q+V_kT_k\tilde A_k^{\top}V_k^{\top}q+V_k\tilde B_k\tilde B_k^{\top}V_k^{\top}q\\
                                       &= V_k(\tilde A_kT_k+T_k\tilde A_k^{\top}+\tilde B_k\tilde B_k^{\top})V_k^{\top}q\\
                                       &= 0
    \end{align*}
    which shows us that the residual is zero, so $V_kT_kV_k^{\top}$ is the exact solution.
  \end{proof}
\end{proposition}
\begin{corollary}
\redii{
From Proposition~\ref{t3} it follows that when $m=1$, the equality $\Range{V_k} = \mathcal{K}_k(A,B)$ holds in Proposition~\ref{t1}.
}
\end{corollary}
\subsection{\green{Restart strategy}}
A problem that occurs in the method described above is that the space $V$ might get quite large.
This means that it can take up a lot of memory, but also that the reduced system, for which we use a dense solver, can become large and take up most of the computation time.
For this reason we implemented a restart strategy, where we reduce the size of $V$ after a certain amount of iterations.
Usually, not all directions that are present in $V$ are equally important, so we just want to keep the most important ones.
We do this by computing the \redii{eigenvectors associated with the largest eigenvalues} of $VTV^{\top}$, which are then used as $V$ in the next iteration of our method.
\blue {
Note that since $V$ is orthonormalized, the nonzero eigenvalues} of $VTV^{\top}$ are the same as the eigenvalues of $T$.
The eigenvectors are given by $VU$, where $U$ are the eigenvectors of $T$, which makes it quite easy to obtain them.

Besides limiting the size of the reduced problem we have to solve, another advantage is that we reduce the rank of the approximate solution, since we only keep the most important components.
This means that we need less memory to store the solution, but also that when we apply the solution, for instance when computing eigenvalues of the solution, we need fewer operations.
To assure \pinkii{that} we have a solution that has a minimal size, we also apply a restart \pinkii{when} \methodname has converged, after which we let it converge once more.
This usually leads to an approximation of lower rank.

A downside of restarting an iterative method is that we lose a lot of convergence properties, like for instance the finite termination property that was shown in Proposition~\ref{t3}.
\redii{Since we keep reducing the size of the search space at the time of a restart}, it might happen that stagnation occurs.
\pinkii{This can be prevented by (automatically) increasing} the tolerance of \pinkii{the vectors that} we retain during a restart.
If we keep doing this repetitively, eventually, the method should still converge.

To implement this restart method, we replaced Line~\ref{alg:stop} in Algorithm~\ref{alg:proj} by Algorithm~\ref{alg:restart}.
\begin{algorithm}[H]
\centering
\parbox{.95\textwidth}{
\begin{framed}
\begin{tabular}{lll}
        \textbf{input:}  & $k$ & Iterations after which to restart.\\
        & $\tau$ & Tolerance for the eigenvalues at a restart.
\end{tabular}
\begin{algorithmic}[1]
  \vspace{1mm}
  \State Set converged to true if the approximated largest eigenvalue is smaller than $\epsilon$
  \If{converged and convergence was already achieved earlier}
  \State \textbf{stop}
  \ElsIf{converged or $\textrm{mod}(j, k)=0$}
  \State Compute the eigenpairs ($\lambda_p$, $q_p$) of $T_j$
  \State $U = [~]$
  \ForAll{eigenvalues $\lambda_p$ larger than $\tau$}
  \State $U = [U, q_p]$
  \EndFor
  \State $V_{j+1} = V_jU$
  \State $\tilde A_{j+1} = U^{\top}\tilde A_j U$
  \State $\tilde M_{j+1} = U^{\top}\tilde M_j U$
  \State $\tilde B_{j+1} = U^{\top}\tilde B_j$
  \EndIf
\end{algorithmic}
\end{framed}}
\caption{Restart method that replaces Line~\ref{alg:stop} Algorithm~\ref{alg:proj}.}
\label{alg:restart}
\end{algorithm}

\section{\green{Application in physical oceanography}}\label{sec:test}
Motivated by the possibility of multiple steady states in the climate system, we have chosen a problem from physical oceanography, involving the Atlantic Ocean circulation.
Its Meridional Overturning Circulation (MOC) is sensitive to freshwater anomalies \cite[]{Stommel1961}.
Freshening of the surface waters in the Nordic and Labrador Seas diminishes the production of deepwater that feeds the deep southward branch of the MOC.
The weakening of the MOC leads to reduced northward salt transport freshening the northern North Atlantic and amplifying the original freshwater perturbation.
A MOC collapse, where a transition to a different steady state (weak MOC) occurs within a few decades, may occur due to the existence of a tipping point (e.g.~due to a saddle-node bifurcation) associated with this salt-advection feedback.
If the MOC is in a multiple steady state regime it may undergo transitions due to the impact of noise.
The stochastic nature of the forcing (wind-stress and/or buoyancy flux) may even lead to transitions \emph{before} the actual saddle-node bifurcation has been reached.
Hence, it is of central importance to determine what perturbations can cause a transition to the weak MOC state, in particular, which spatial and temporal correlations in the noise are most effective.

We use the spatially quasi two-dimensional model of the Atlantic MOC as described in \cite{Toom2011}.
In the model, there are two active tracers: temperature $T$ and salinity $S$, which are related to the density $\rho$ by a linear equation of state
\begin{equation}
\rho =\rho_0\brx{1-\alpha_T\brx{T -T_0}+\alpha_S\brx{S -S_0}},
\label{eq:eos}
\end{equation}
where $\alpha_T$ and $\alpha_S$ are the thermal expansion and haline contraction coefficients, respectively, and $\rho_0$, $T_0$, and $S_0$ are reference quantities.
The numerical values of the fixed model parameters are summarized in \tabref{modelpars}.
\begin{table*}[tph]
\centering
  \begin{tabular}{|rclcll|rclcll|}
\hline
  $D           $   &=& 4.0     &$\cdot$& $10^3$         & m                            & $H_m       $   &=& 2.5    & $\cdot$& $10^2$       & m                 \\
  $r_0        $  &=& 6.371 &$\cdot$& $10^6$         & m                             &$T_0         $   &=& 15.0 &              &                  & $^\circ$C       \\
  $g           $   &=& 9.8      &             &                       & m s$^{-2}$           &$S_0         $   &=& 35.0&               &                  & psu                   \\          
  $A_H      $  &=& 2.2         &$\cdot$& $10^{12}$ & m$^2$s$^{-1}$    & $\alpha_T$   &=& 1.0    &$\cdot$& $10^{-4}$&K$^{-1}$         \\
  $A_V      $  &=& 1.0      &$\cdot$& $10^{-3}$     & m$^2$s$^{-1}$    & $\alpha_S$   &=& 7.6    &$\cdot$& $10^{-4}$&psu$^{-1}$      \\
  $K_H      $  &=& 1.0     &$\cdot$& $10^3$          & m$^2$s$^{-1}$    & $\rho_0     $  &=& 1.0    &$\cdot$& $10^3$    & kg m$^{-3}$   \\
  $K_V      $  &=& 1.0      &$\cdot$& $10^{-4}$     & m$^2$s$^{-1}$    & $\tau          $   &=& 75.0 &             &                  &days                 
  \\\hline
  \end{tabular}
\caption{Fixed model parameters of the two-dimensional ocean model.}
\label{tab:modelpars}
\end{table*}

In order to eliminate longitudinal dependence from the problem, we consider a purely buoyancy-driven flow on a non-rotating Earth.
We furthermore assume that inertia can be neglected in the meridional momentum equation.
The mixing of momentum and tracers due to eddies is parameterized by simple anisotropic diffusion.
In this case, the zonal velocity as well as the longitudinal derivatives are zero and the equations for the meridional velocity $v$, vertical velocity $w$, pressure $p$, and the tracers $T$ and $S$ are given by
\bse
0=-\frac{1}{\rho_0 r_0}\pd{p}{\theta}+A_V\pd{^2v}{z^2}+
\frac{A_H}{r_0^2}\brx{\frac{1}{\cos\theta}\pd{}{\theta}\brx{\cos\theta\pd{v}{\theta}}+\brx{1-\tan^2\theta}v}, \\
0=-\frac{1}{\rho_0}\pd{p}{z}+g\brx{\alpha_T T-\alpha_S S}, \\
0=\frac{1}{r_0\cos\theta}\pd{v\cos\theta}{\theta}+\pd{w}{z}, \\
\pd{T}{t}+\frac{v}{r_0}\pd{T}{\theta}+w\pd{T}{z}=
\frac{K_H}{r_0^2\cos\theta}\pd{}{\theta}\brx{\cos\theta\pd{T}{\theta}}+K_V\pd{^2T}{z^2}+\mathrm{CA}(T), \\
\pd{S}{t}+\frac{v}{r_0}\pd{S}{\theta}+w\pd{S}{z}=
\frac{K_H}{r_0^2\cos\theta}\pd{}{\theta}\brx{\cos\theta\pd{S}{\theta}}+K_V\pd{^2S}{z^2}+\mathrm{CA}(S). 
\label{eq:2Deqs}
\ese
Here $t$ is time, $\theta$ latitude, $z$ the vertical coordinate, $r_0$ the radius of Earth, $g$ the acceleration due to gravity, $A_H$ ($A_V$) the horizontal (vertical) eddy viscosity, and $K_H$ ($K_V$) the horizontal (vertical) eddy diffusivity.
The terms with CA represent convective adjustment contributions.

The equations are solved on an equatorially symmetric, flat-bottomed domain.
The basin is bounded by latitudes $\theta = - \theta_N$ and $\theta=\theta_N = 60^\circ$ and has depth $D$.
Free-slip conditions apply at the lateral walls and at the bottom.
Rigid lid conditions are assumed at the surface and atmospheric pressure is neglected.
The wind stress is zero everywhere, and ``mixed'' boundary conditions apply for temperature and salinity,
\begin{subequations}
\begin{align}
K_V\pd{T}{z}&=\frac{H_m}{\tau}\brx{\bar{T}(\theta)-T} , \\
K_V\pd{S}{z}&= S_0 F_s(\theta).                        
\end{align}
\label{eq:2Dforcing}
\end{subequations}
This formulation implies that temperatures in the upper model layer (of depth $H_m$) are relaxed to a prescribed temperature profile $\bar{T}$ at a rate $\tau^{-1}$, while salinity is forced by a net freshwater flux $F_s$, which is converted to an equivalent virtual salinity flux by multiplication with $S_0$.
The specification of the CA terms is given in \cite{Toom2011}.

In \secref{bifdiag} below, we first determine the bifurcation diagram of the deterministic model using pseudo-arclength continuation methods.
In the next sections, the case with stochastic freshwater forcing is considered, focusing on validation of the new methods (\secref{validation}), comparison with other methods (\secref{comp}), numerical aspects (\secref{scal}) and potential extensions (\secref{3d}).

\subsection{Bifurcation diagram} \label{sec:bifdiag}

For the deterministic case, we fix the equatorially symmetric surface forcing as 
\bse
\bar{T}(\theta) &=& 10.0\cos\brx{\pi\theta/\theta_N},\\
F_s(\theta) = {\bar F}_s(\theta) &=&  \mu F_0 \frac{\cos(\pi\theta/\theta_n)}{\cos(\theta)}, 
\label{e:FSd}
\ese  
where $\mu$ is the strength of the mean freshwater forcing (which we take as bifurcation parameter) and $F_0 = 0.1$ m yr$^{-1}$ is a reference freshwater flux.

The equations are discretized on a latitude-depth equidistant $n_x \times n_y \times n_z$ grid using a second-order conservative central difference scheme.
An integral condition expressing the overall conservation of salt is also imposed, as the salinity equation is only determined up to an additive constant.
The total number of degrees of freedom is $n = 6 n_x n_y n_z$, as there are six unknowns per point.
The standard spatial resolution used is $n_x = 4, n_y= 32, n_z = 16$ and the solutions are uniform in the zonal direction, with the zonal velocity $u = 0$.

\begin{figure}[ht]
  \begin{subfigure}[h]{0.48\textwidth}
    \centering\includegraphics[width=\textwidth]{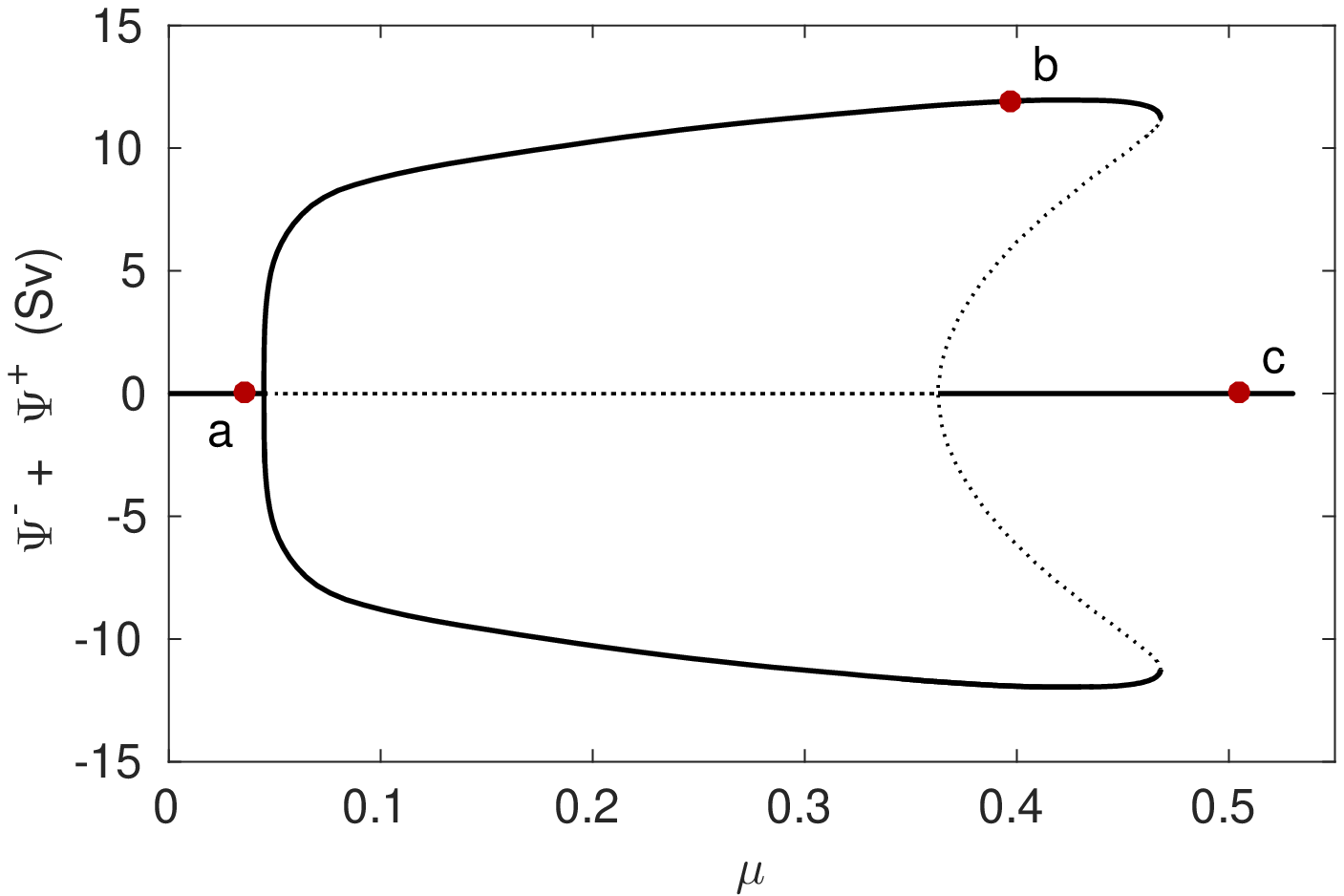}
    \caption{Bifurcation diagram where a solid line is stable and dashed line is unstable.}
  \end{subfigure}
  \hspace{0.04\textwidth}
  \begin{subfigure}[h]{0.48\textwidth}
    \centering\includegraphics[width=\textwidth]{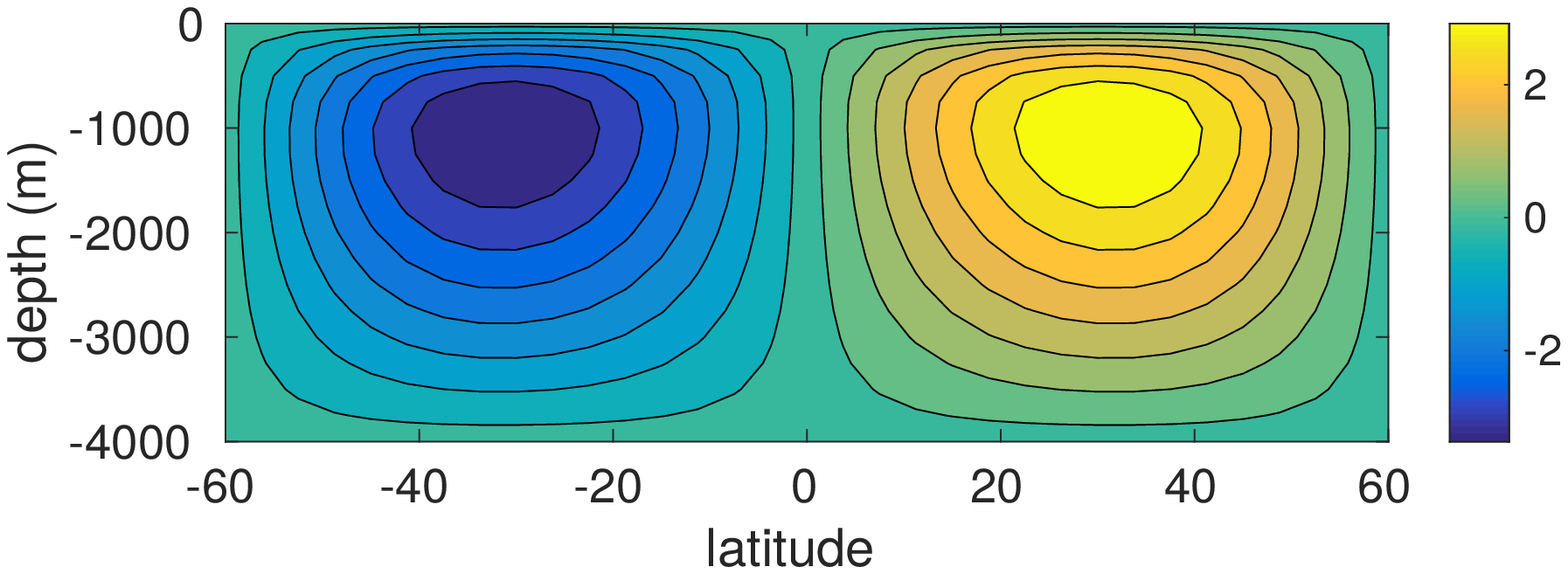}
    \centering\includegraphics[width=\textwidth]{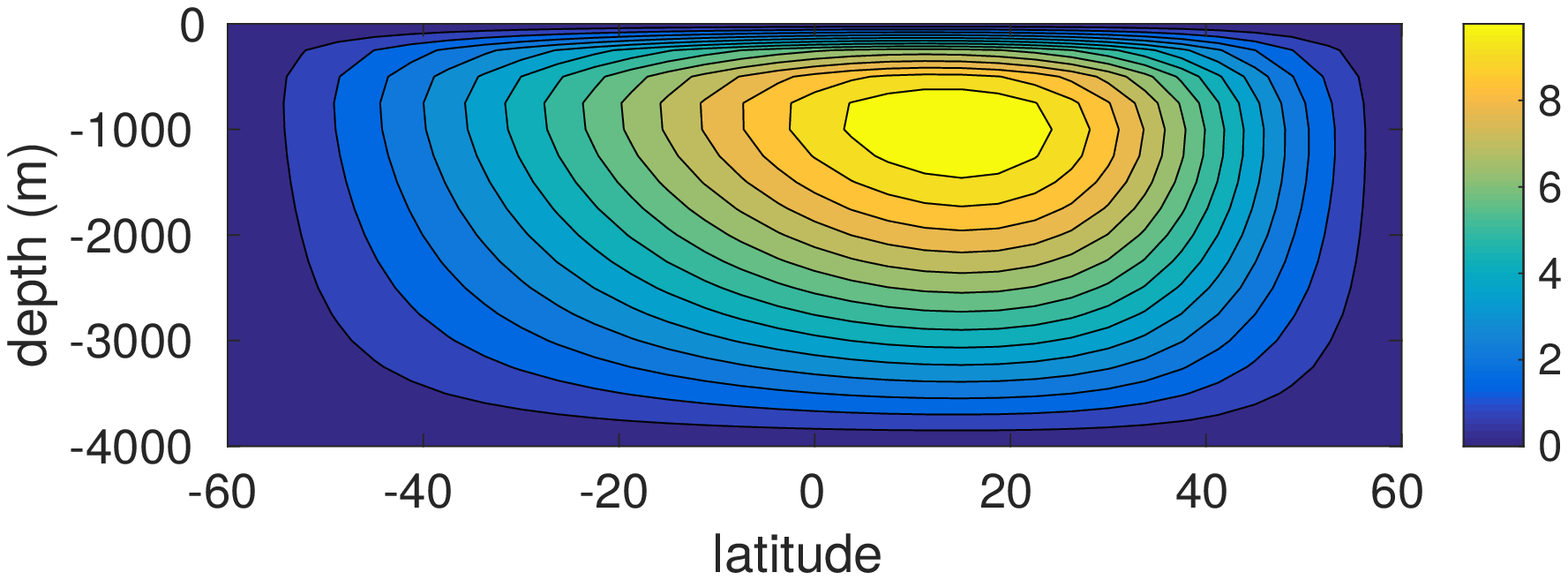}
    \centering\includegraphics[width=\textwidth]{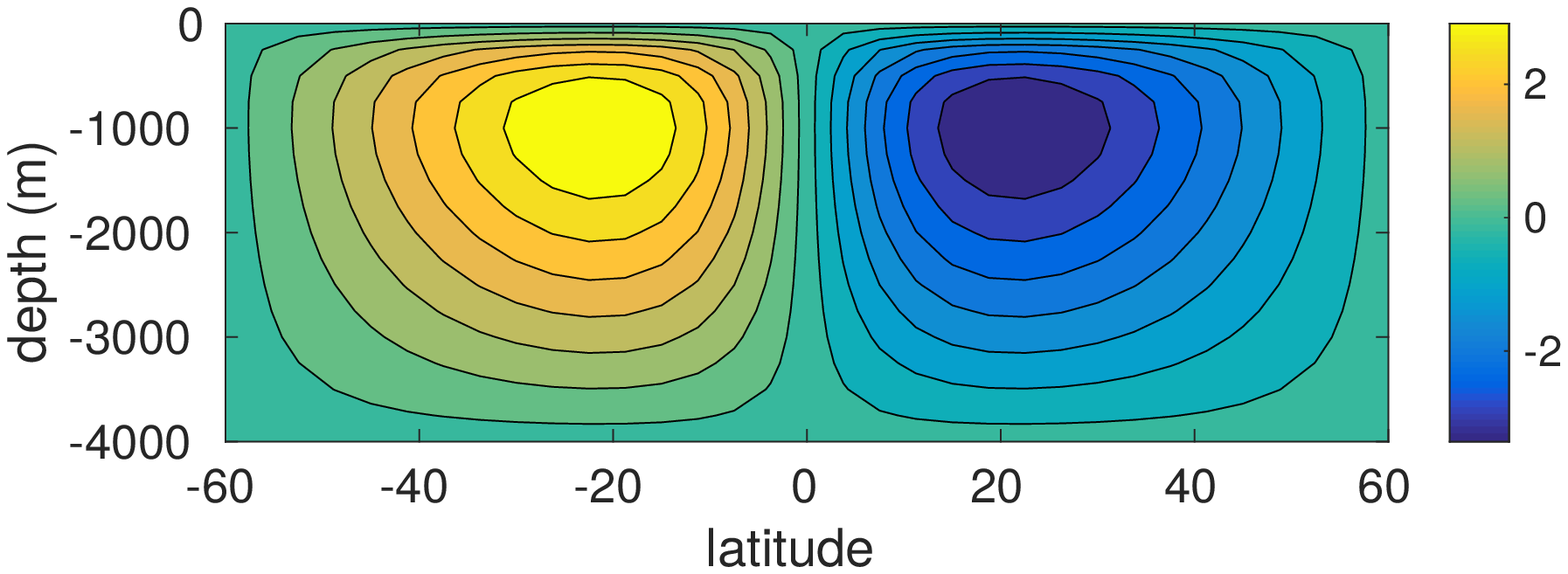}
    \caption{The streamfunction at $\mu_a$, $\mu_b$ and $\mu_c$ respectively.}
  \end{subfigure}
  \caption{(a) Bifurcation diagram of the deterministic equatorially symmetric 2D MOC 
  model, with the forcing as in Eq.~(\ref{e:FSd}). (b) Streamfunction pattern at $\mu_a$, $\mu_b$ and $\mu_c$ respectively.
}
\label{fig:bifdia}
\end{figure}

The bifurcation diagram of the deterministic model for parameters as in Tab.~1 is shown in \figref{bifdia}a.
On the $y$-axis, the sum of the maximum ($\Psi^+$) and minimum ($\Psi^-$) values of the meridional streamfunction $\Psi$ is plotted, where $\Psi$ is defined through
\be
\pd{\Psi}{z} = v  \cos\theta, ~ -  \frac{1}{r_0\cos\theta} \pd{\Psi}{\theta} = w .
\ee
For the calculation of the transports, the basin is assumed to have a zonal width of $64^\circ$.
The value of $\Psi^+ + \Psi^-$ is zero when the MOC is symmetric with respect to the equator.

For small values of $\mu$, a unique equatorially anti-symmetric steady MOC exists of which a pattern at location $a$ is shown in \figref{bifdia}b.
This pattern destabilizes at a supercritical pitchfork bifurcation and two asymmetric pole-to-pole solutions appear.
An example of the MOC at location $b$ in \figref{bifdia}b shows a stronger asymmetric overturning with sinking in the northern part of the basin.
The pole-to-pole solutions cease to exist beyond a saddle-node bifurcation near $\mu = 0.47$ and both branches connect again with the anti-symmetric solution at a second supercritical pitchfork bifurcation.
At this bifurcation, the anti-symmetric solution with equatorial sinking (see MOC at location $c$ in \figref{bifdia}b) appears which is stable for larger values of $\mu$.
The value of $\mu$ at the point $b$, $\mu_b=0.40$, will be our reference freshwater forcing.

\subsection{Stochastic freshwater forcing} \label{sec:forcing}
\green {
\pinkii{The freshwater forcing is chosen as}
\bae\label{FS}
F_s(\theta, t) = (1+ \sigma ~\zeta(\theta, t)) \bar{F}_s(\theta) 
\eae
where $\zeta(\theta, t)$ represents zero-mean white noise with a unit standard deviation, i.e., with $\E[\zeta(\theta, t)] = 0$, $\E[\zeta(\theta, t) \zeta(\theta, s)] = \delta(\theta, t-s)$ and $\E[\zeta(\theta, t) \zeta(\eta, t)] = \delta(\theta-\eta, t)$.
The constant $\sigma$ is the standard deviation of the noise which we set to $\sigma=0.1$.
}
In this case, the noise matrix $B$ in Eq.~\eqref{SDAElin} simply represents additive noise which is $(i)$ only active in the freshwater component, $(ii)$ only present at the surface, $(iii)$ meridionally uncorrelated (unless stated otherwise), and $(iv)$ has  magnitude $\sigma$ of $10\%$ of the deterministic freshwater forcing amplitude at each latitude $\theta$ (see Eq.~\eqref{FS}).

\section{\green{Results}}\label{sec:results}
Using the available Jacobian $A$ of the deterministic continuation, the mass matrix $M$, which is a diagonal matrix with non-zero elements in the $T$ and $S$ rows, and the forcing $B$ as described above, we can determine the local probability distribution of a steady state using the generalized Lyapunov equation \eqref{schurlyapunov}.
We use grid sizes of $4 \times n_y \times 16$, where $n_y$ is a varying amount of grid points in the meridional direction, 16 is the amount of grid points in the vertical direction and \pinkii{there are} 4 grid points in the zonal direction.
Since our model is two-dimensional, both the forcing and the solution will be constant in this direction.
For the forcing, this means that $B$ contains $n_y$ vectors with the forcing as described in Eq.~\eqref{FS}.

For \methodname, we use the algorithm as described in \secref{methods} and always expand with $m = 3$ vectors per iteration unless stated otherwise.
When comparing to other Lyapunov solvers, which were mostly written in Matlab, we use a Matlab implementation of \methodname (\secref{comp}), but when we solve larger systems, we prefer to use a C++ implementation (all other sections).
Computations are performed on one node of Peregrine, the HPC cluster of the University of Groningen.
We use this machine to be able to make fair comparisons with other methods, which use large amounts of memory.
Peregrine has nodes with 2 Intel Xeon E5 2680v3 CPUs (24 cores at 2.5GHz) and each node has 128 GB of memory.
Only one core is used in the results below to be able to make fair comparisons.

\subsection{\green{Comparison with stochastically forced time forward simulation}} \label{sec:validation}
\begin{figure}[th!]
  \begin{subfigure}[h]{\textwidth}
    \centering\includegraphics[width=\textwidth]{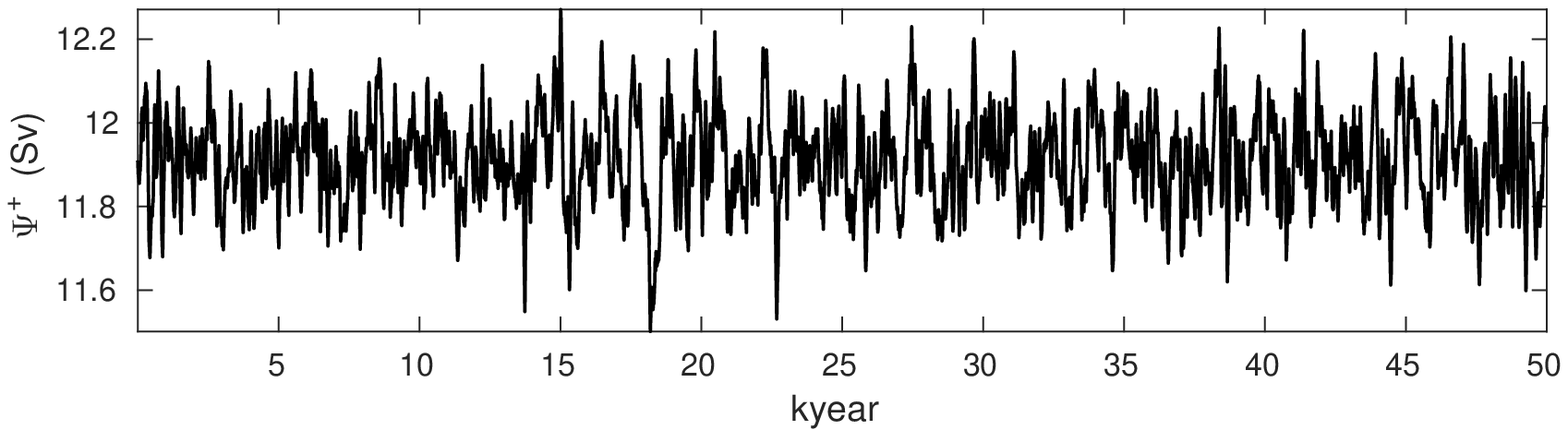}
    \caption{Time series of the MOC maximum $\Psi^+$.}
    \label{fig:time1}
  \end{subfigure}
  \begin{subfigure}[h]{0.48\textwidth}
    \centering\includegraphics[width=\textwidth]{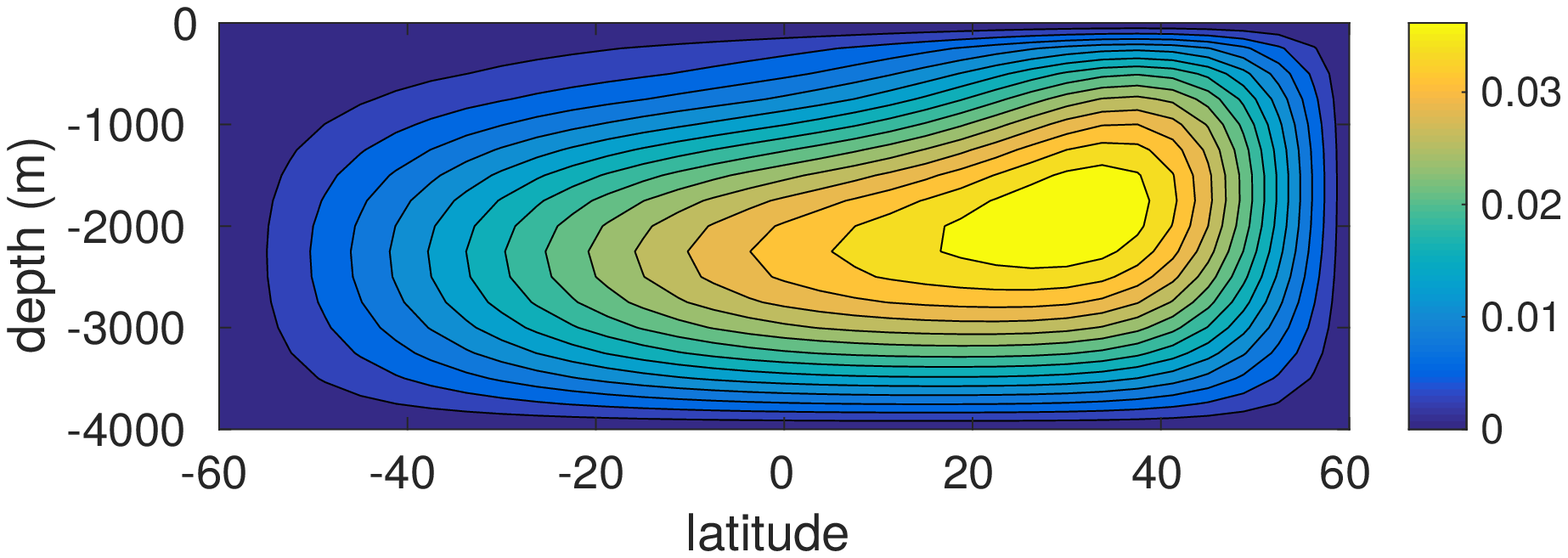}
    \centering\includegraphics[width=\textwidth]{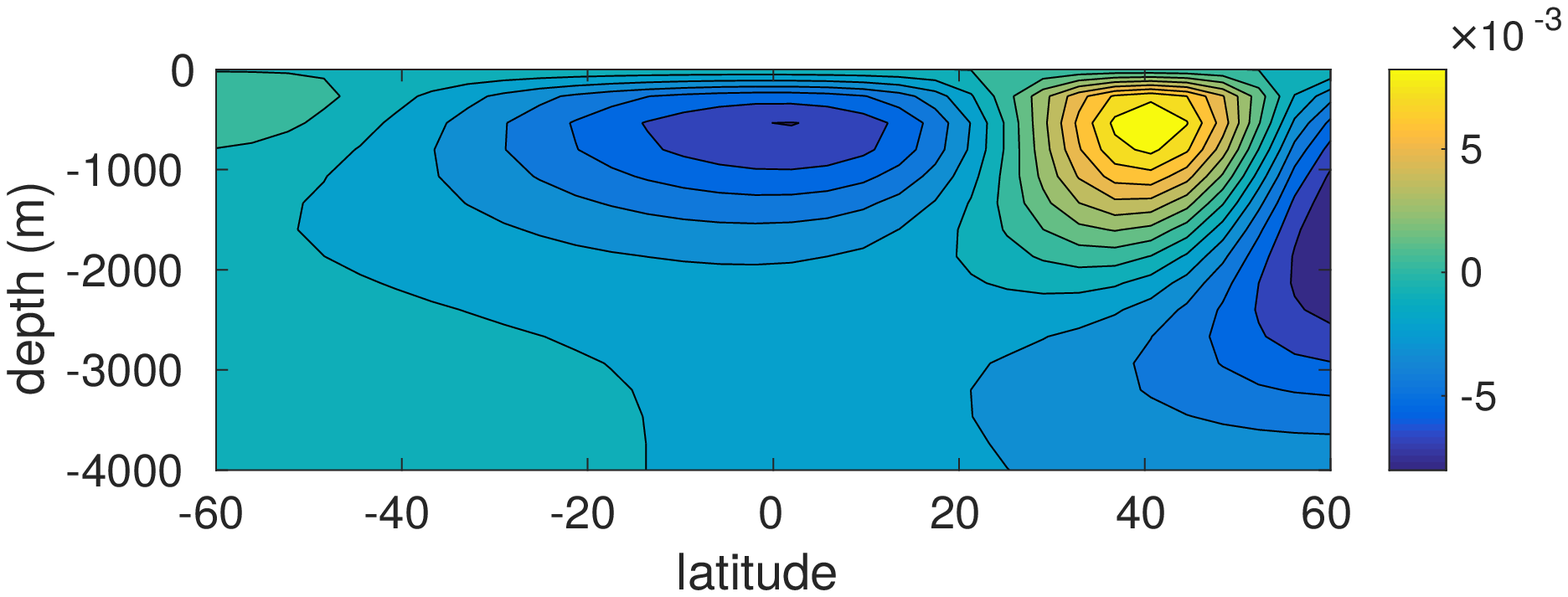}
    \centering\includegraphics[width=\textwidth]{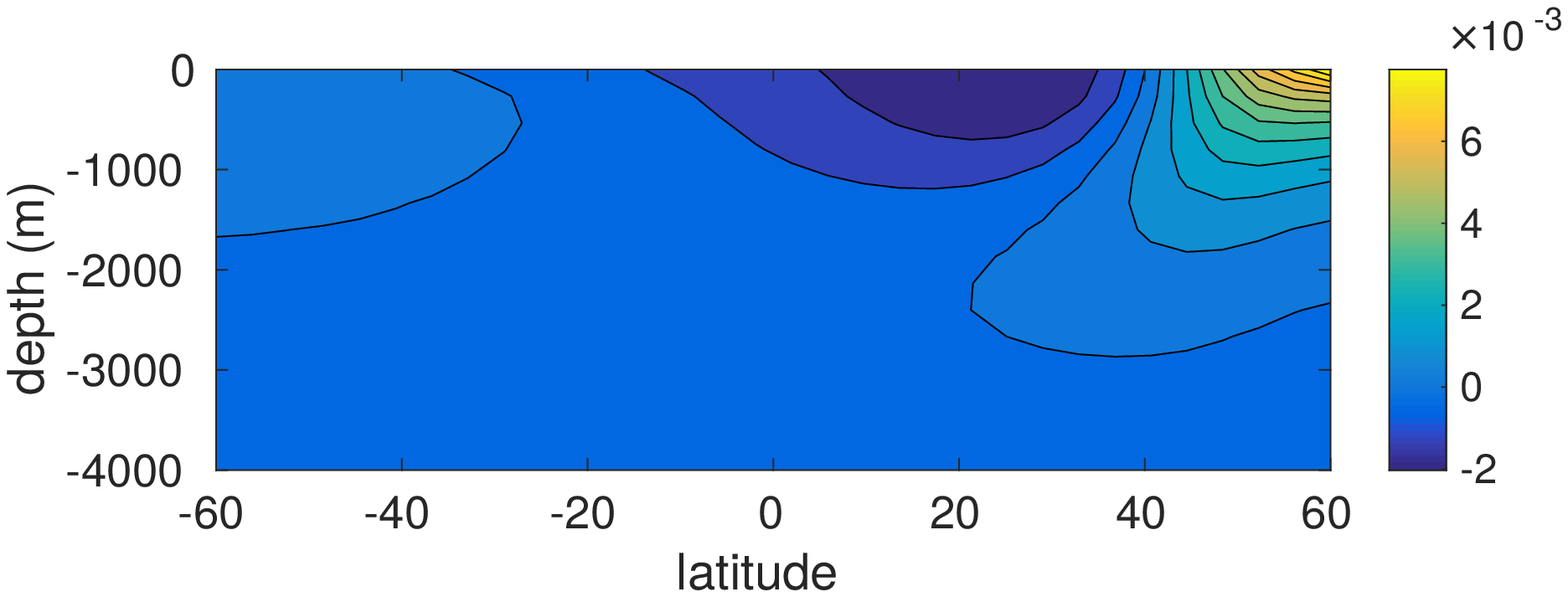}
    \caption{The first EOF
    }
    \label{fig:timeeofs1}
  \end{subfigure}
  \hspace{0.04\textwidth}
  \begin{subfigure}[h]{0.48\textwidth}
    \centering\includegraphics[width=\textwidth]{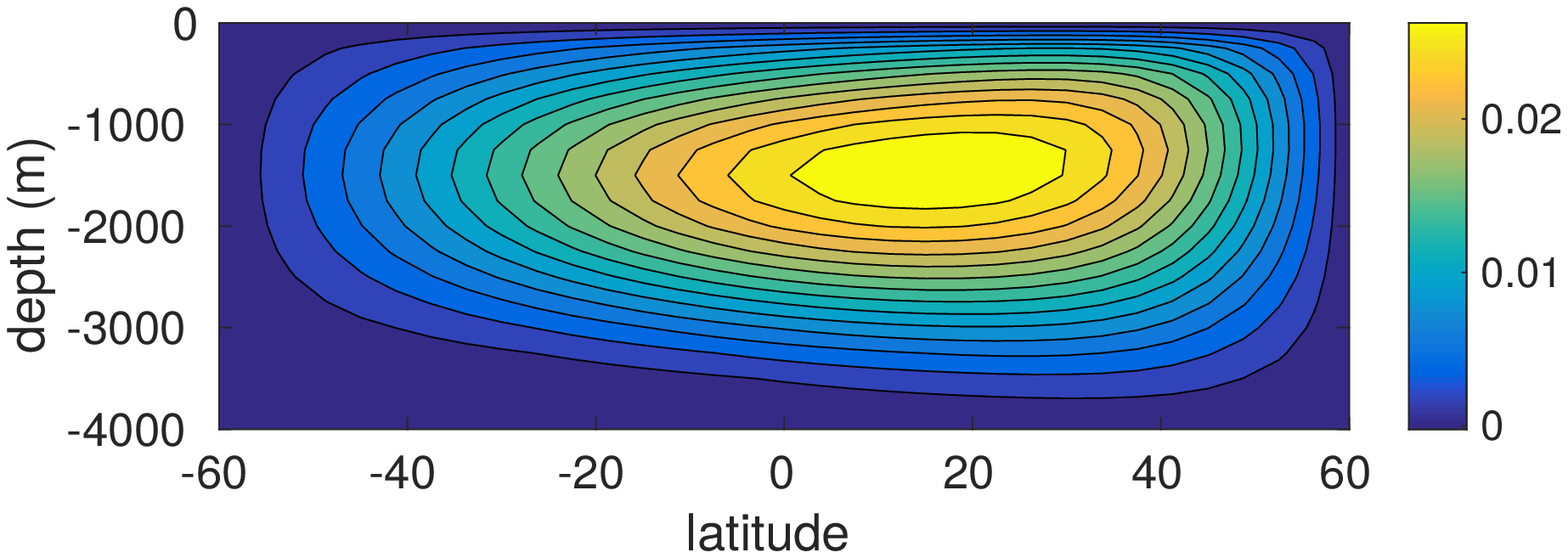}
    \centering\includegraphics[width=\textwidth]{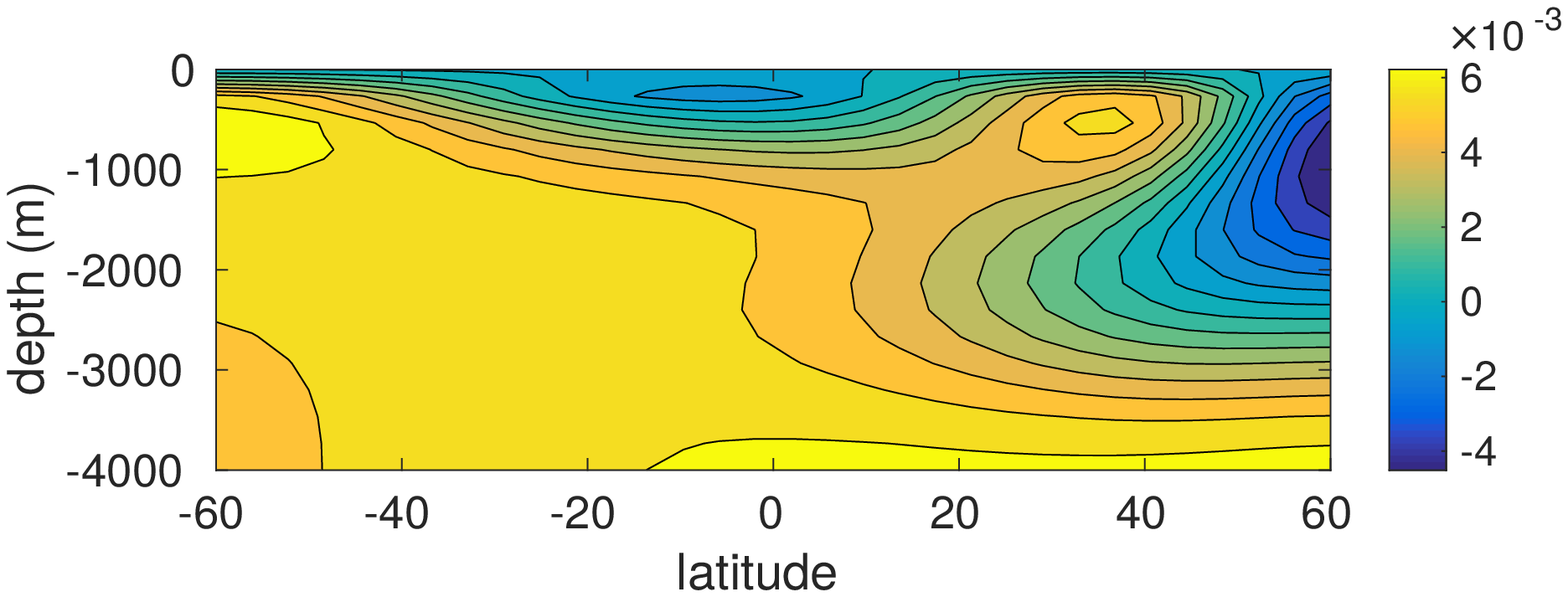}
    \centering\includegraphics[width=\textwidth]{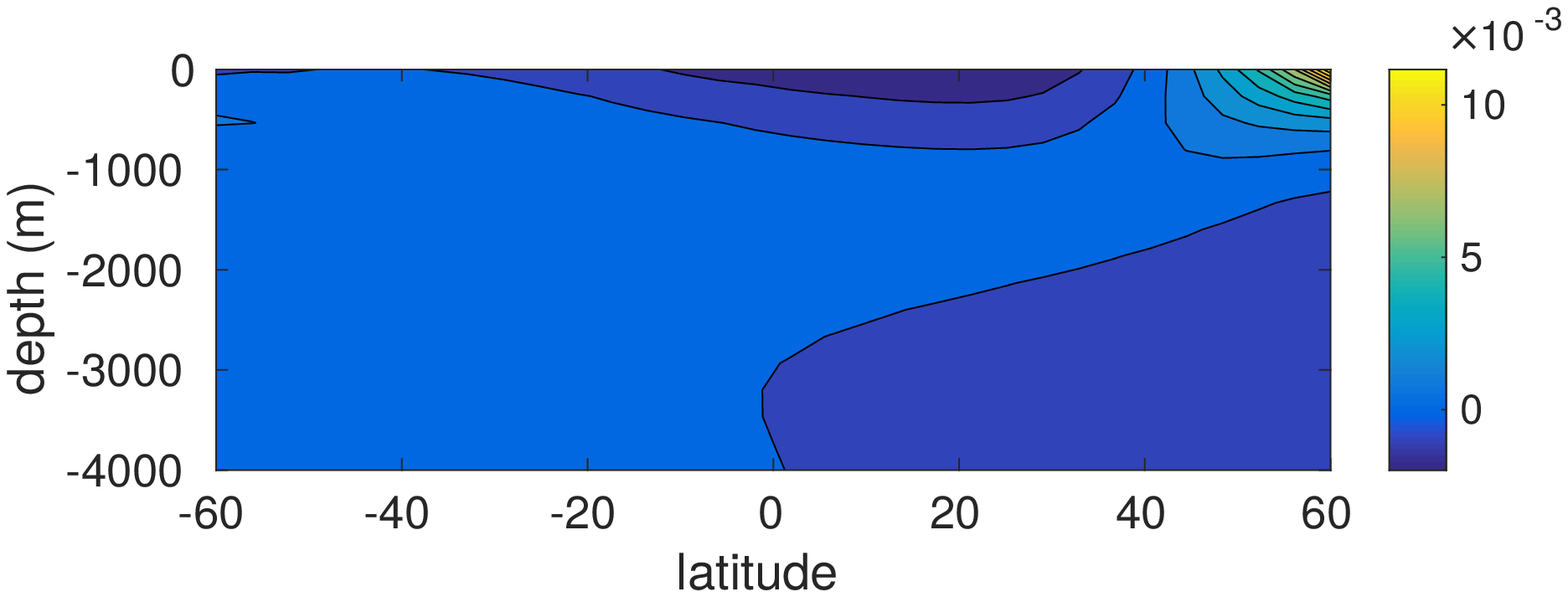}
    \caption{The second EOF.
    }
    \label{fig:timeeofs2}
  \end{subfigure}
  \caption{(a) Time series of the maximum MOC strength $\Psi^+$ for a 50,000 years simulation of the model for $\mu = \mu_b$ under the freshwater forcing as in Eq.~\eqref{FS}.
The variability has a slight negative skewness of $-0.05$ and a kurtosis of $3.04$.
(b) Patterns of $\Psi$ (top), isothermals (middle) and isohalines (bottom) for the first EOF of this simulation.
(c) Similar to (b) but for the second EOF.
\label{fig:timeeofs}}
\end{figure}

A first check of the correctness of the approximate solution of the generalized Lyapunov equations is obtained by comparing the empirical orthogonal functions (EOFs) and weighted eigenvalues of the covariance matrix \pinkii{that} we get from both the Lyapunov solver and a stochastically forced time forward simulation at $\mu = \mu_b$, similar to those performed in \cite{Mheen2013}.
This time series (for $n_y=32$) is plotted in \figref{time1} and shows that $\Psi^+$ fluctuates around the mean MOC value at $\mu_b$.
The patterns of the MOC, the temperature field and the salinity field of both EOF1 and EOF2 are shown in \figref{timeeofs}b and c.

\begin{figure}[ht]
  \begin{subfigure}[h]{0.48\textwidth}
    \centering\includegraphics[width=\textwidth]{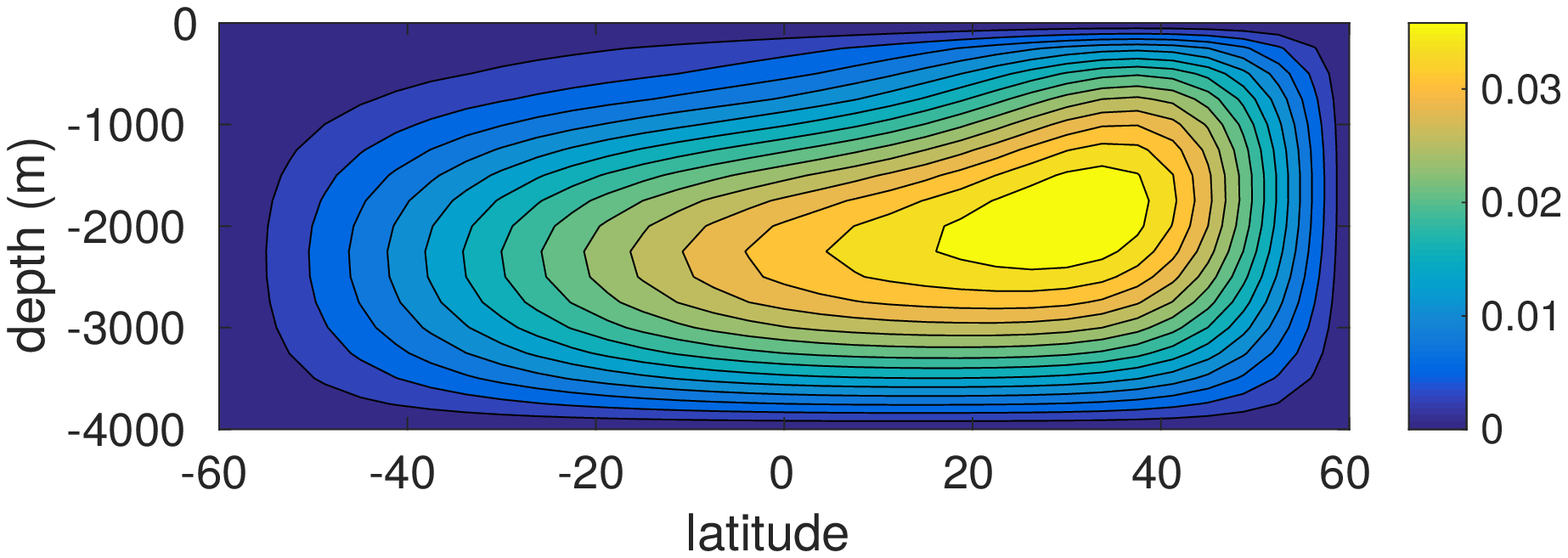}
    \centering\includegraphics[width=\textwidth]{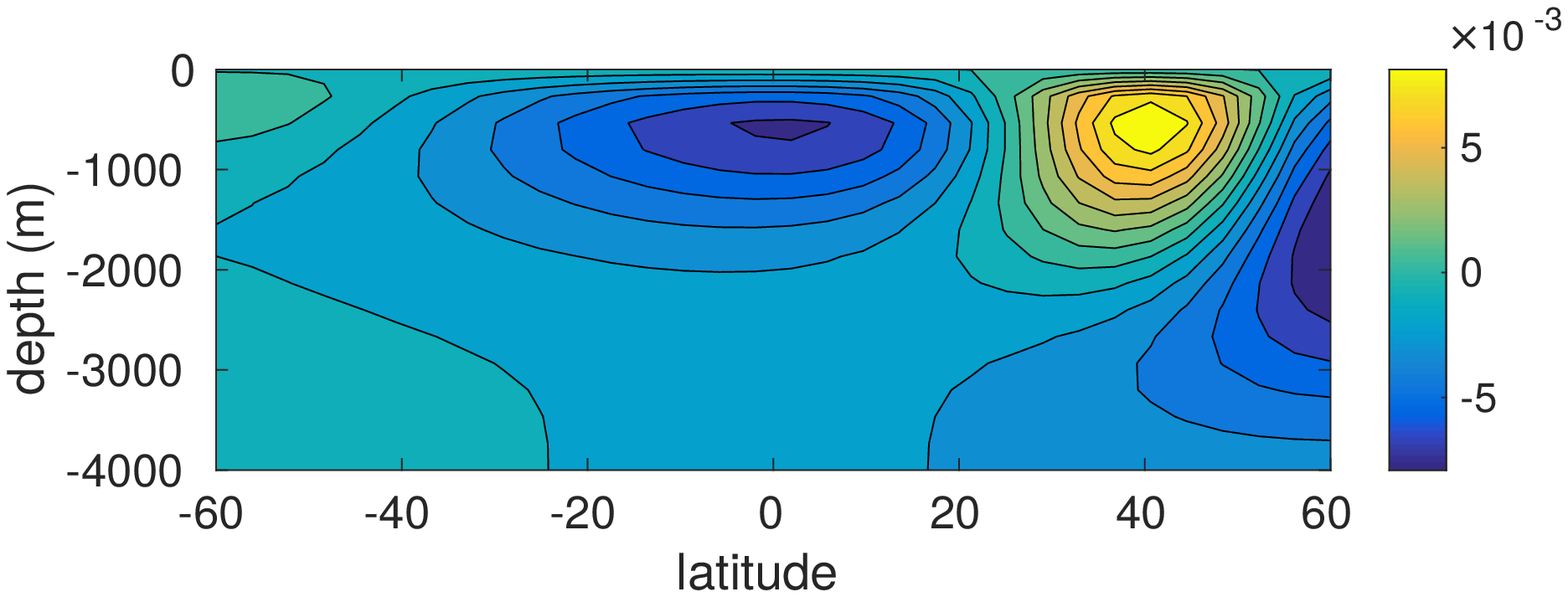}
    \centering\includegraphics[width=\textwidth]{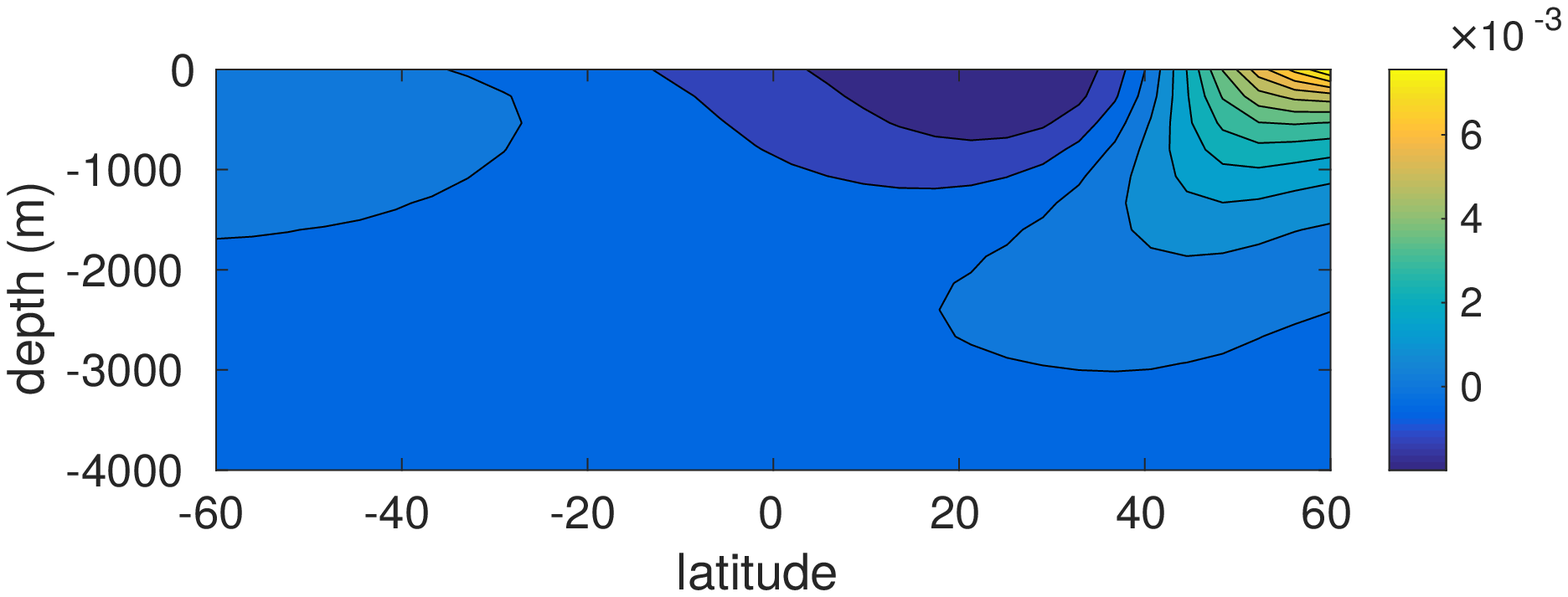}
    \caption{The first EOF
    }
    \label{fig:lyapeofs1}
  \end{subfigure}
  \hspace{0.04\textwidth}
  \begin{subfigure}[h]{0.48\textwidth}
    \centering\includegraphics[width=\textwidth]{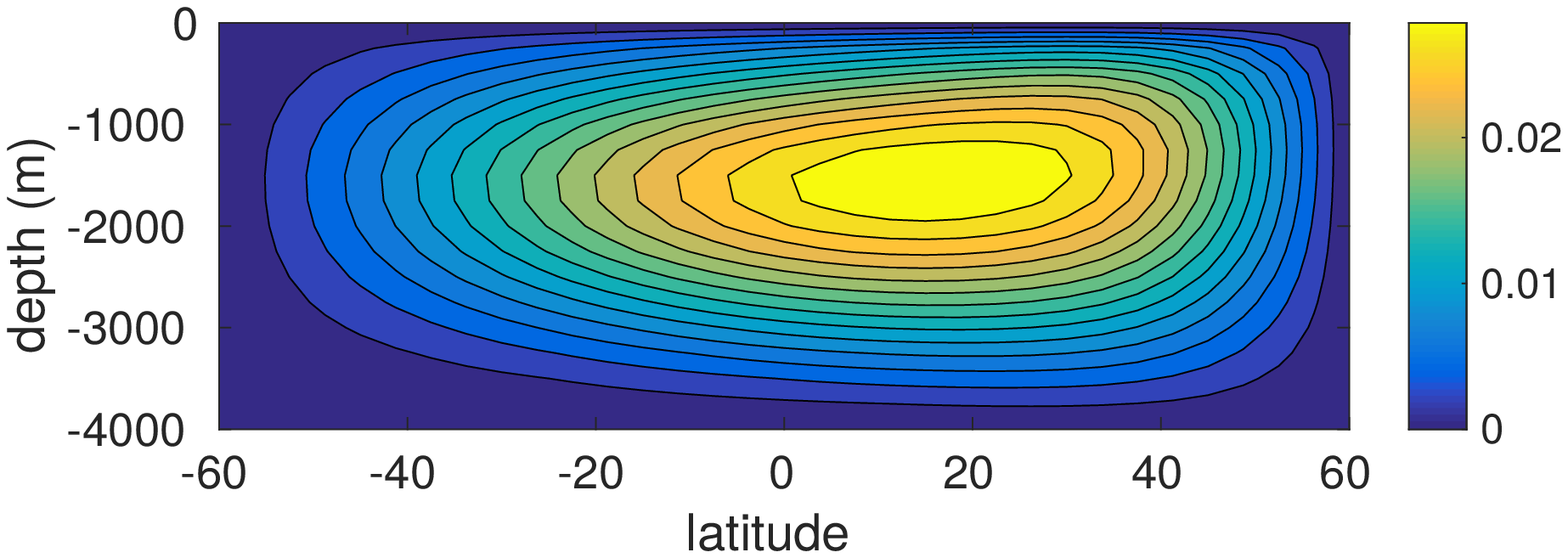}
    \centering\includegraphics[width=\textwidth]{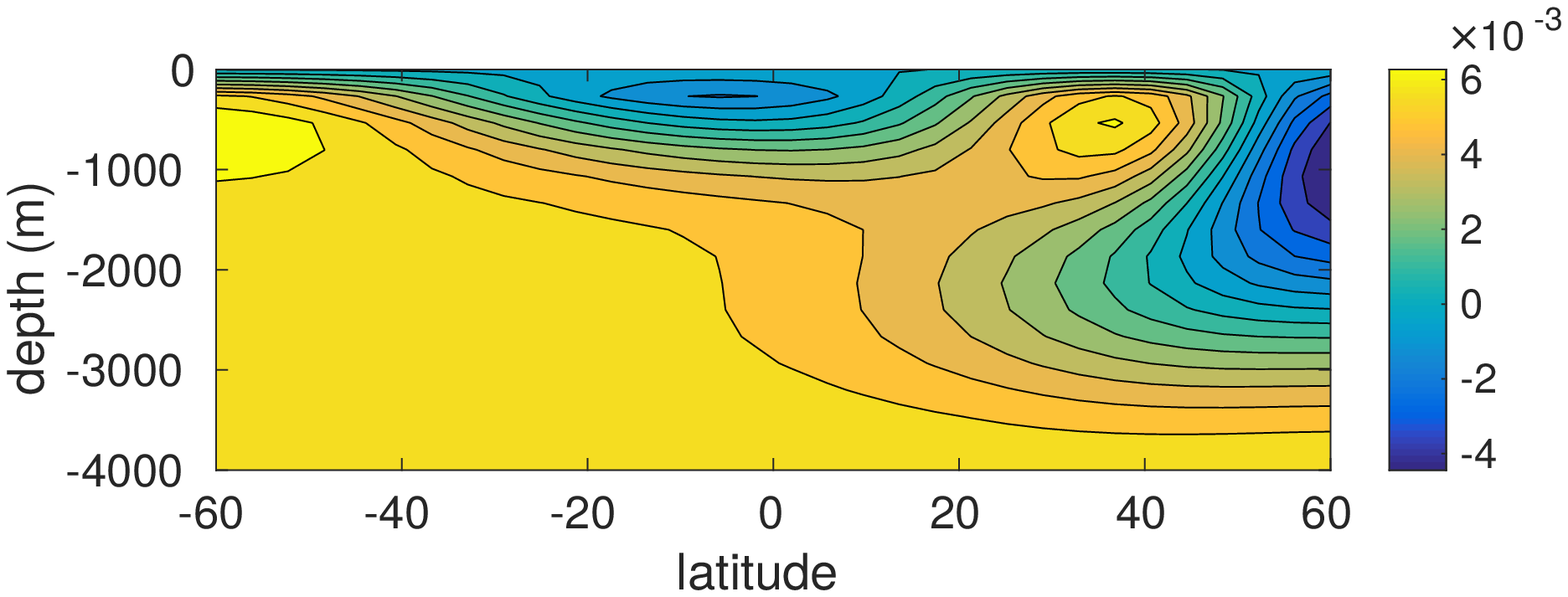}
    \centering\includegraphics[width=\textwidth]{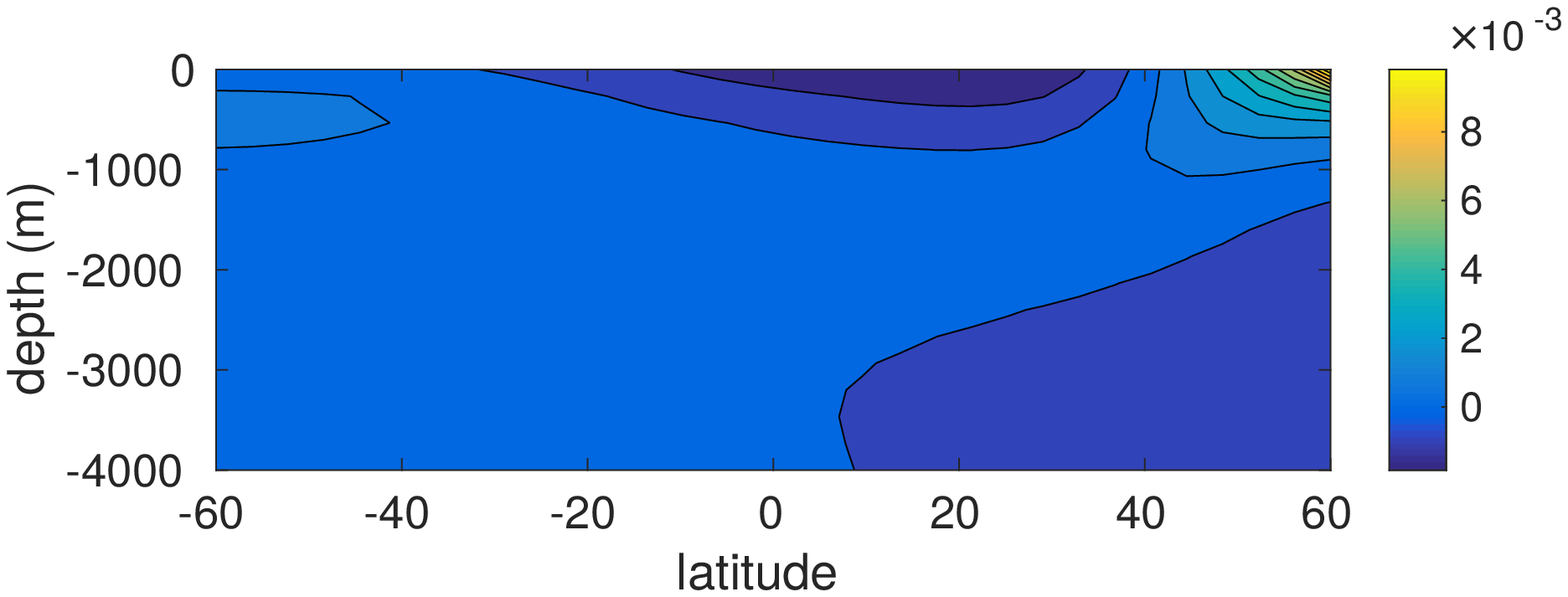}
    \caption{The second EOF.
    }
    \label{fig:lyapeofs2}
  \end{subfigure}
  \caption{(a) Patterns of $\Psi$ (top),   isothermals (middle) and isohalines (bottom) respectively for the first EOF 
 obtained by solving the generalized Lyapunov equation. (b) Similar to  (a) but for the second EOF.  
 \label{fig:lyapeofs}}
\end{figure}

The eigensolutions of the generalized Lyapunov equation, also for $n_y=32$, are shown in \figref{lyapeofs}.
Comparing \figref{timeeofs} with \figref{lyapeofs}, we see that the results from the transient flow computation and the approximate solution of the generalized Lyapunov equation look very similar.
Also, the eigenvalues we find with both methods are very similar, as can be seen in \tabref{eigs}.
\pinkii{The fact that they are not exactly the same is most likely due the fact that the time series takes a long time to converge to a statistical steady state.}
Since the eigenvalues are really close nevertheless, we can indeed use \methodname to compute an estimate of the local probability distribution of steady states of the MOC.

As another check, we use the Bartels-Stewart algorithm \cite{Bartels1972}, which is a dense solver of which the solution time increases with $\mathcal{O}(n^3)$ and the required memory with $\mathcal{O}(n^2)$.
The implementation we use is sb03md from the SLICOT library \cite{Benner1999}.
Results from this method (\tabref{eigs}, also for $n_y=32$) confirm that our solution method provides correct solutions.

\begin{table}[ht]
\centering
\begin{tabular}{l|llll}
Method & $\lambda_1$ & $\lambda_2$ & $\lambda_3$ & $\lambda_4$\\\hline
\methodname & 0.677 & 0.176 & 0.078 & 0.033\\
Dense Lyapunov & 0.677 & 0.176 & 0.079 & 0.033\\
Time series & 0.679 & 0.170 & 0.082 & 0.033
\end{tabular}
\caption{First four weighted eigenvalues of the covariance matrix.
For the \methodname solver $\Vert R \Vert_2 / \Vert B B^T\Vert_2 < 10^{-2}$ was used as stopping criterion.
\label{tab:eigs}}
\end{table}
\subsection{Comparison with other Lyapunov solvers}\label{sec:comp}\llabel{comp}
\blue{
In this section, we compare the results of \methodname to those obtained with standard implementations of the Extended Krylov (EKSM) \cite{Simoncini2007}, Projected Extended Krylov (PEKSM) \cite{Stykel2012}, \blue{Rational Krylov (RKSM) \cite{Druskin2011}, Tangential Rational Krylov (TKRSM) \cite{Druskin2014}} and Low-rank ADI (LR-ADI) \cite{Penzl1999} methods.
For the LR-ADI method, results with the heuristic shifts from \cite{Penzl1999} and self-generating shifts based on a Galerkin projection from \cite{Benner2014} are reported.
Implementations of the Extended Krylov and Rational Krylov methods were obtained from the website of Simoncini \cite{simonciniweb}, whereas the LR-ADI method from M.E.S.S. was used \cite{Saak2016}.
}

\green{
What we want to show is that \methodname can be competitive without requiring linear system solves in every iteration, which all the methods we compare to require.
However, convergence properties of the other methods might be better since they use these linear system solves.
For this reason, we also add experiments with our method, where we expand the search space by $S^{-1}r_i$, where $r_i$ are the eigenvectors of the residual \redii{associated with} the largest eigenvalues.
From now on we will refer to this method as Inverse \methodname.
To be able to better compare to Extended Krylov, we could also expand our search space by $[r_i, S^{-1}r_i]$, but some preliminary experiments showed that Inverse \methodname performs slightly better.
}

\blue{
For Inverse \methodname, Extended/Rational Krylov and ADI based methods, the linear system solves of the form  $Sy=x$ are implemented by first computing an LU-factorization of $A$ using \llabel{umfpack}UMFPACK whenever possible (available from \texttt{lu} in Matlab), and then solving the system 
\begin{align*}
  \begin{pmatrix}
    A_{11} & A_{12}\\A_{21} & A_{22}
  \end{pmatrix}
  \begin{pmatrix}
    \tilde y\\y
  \end{pmatrix}
  =
  \begin{pmatrix}
    0\\x
  \end{pmatrix}.
\end{align*}
This is similar to what has been used in LR-ADI methods in \cite{Freitas2008}.
Alternatives would be first computing $S$ and then making a factorization of $S$, which is not feasible since $S$ tends to be quite dense, or using an iterative method where we need repeated applications of $S$, which includes solving a system with $A_{11}$.
Both alternatives tend to be slower than the method described above if we add up factorization and solution times.
}

\blue {
For the (Tangential) Rational Krylov methods we precompute the initial values $s_0^{(1)}$ and $s_0^{(2)}$ as the eigenvalues of $S$ with the smallest and largest real part.
The time it required to compute the eigenvalues is not included in the results.
}

\blue{
\llabel{proj}For the Projected Extended Krylov method, we followed \cite{Marz1996} for obtaining the spectral projectors that are required and implemented the method according to \cite{Stykel2012}.
The advantage of this method is that it does not require the Schur complement as we described above.
The disadvantage is that instead spectral projectors are required.
For this method, we have to compute the projected matrix $V^{\top}AV$ explicitly, and not implicitly as suggested in \cite{Stykel2012} to obtain sufficient accuracy.
Without this, the Lyapunov solver that was used for the small projected Lyapunov equation would fail to find a solution.
For comparison, we implemented the same projection method that was used in \cite{Stykel2012} also in \methodname, where we chose to expand the basis with $A^{-1}r_i$.
}

\pink{
\llabel{restol}For all methods we use the same stopping criterion, namely that \pinkii{we require} the relative residual $\rho = \Vert R \Vert_2 / \Vert BB^{\top} \Vert_2 < \epsilon$, where $\epsilon = 10^{-2}$, which is sufficient for our application.
The resulting absolute residual norm will be around $10^{-5}$ for the MOC problem with $n_y=32,64,128$.
We also show the final relative residual in our results.}
\red {
\llabel{init}For \methodname we use a random $V_1$ as initial guess, since this seemed to give slightly better results than taking $B$ as initial guess.
Using a random initial guess also shows that even if storing $B$ in a dense way requires too much memory, we are still able to start our method even if the other methods cannot.
During a continuation run, we can of course do much better than a random initial guess by using the approximate solution space from the previous continuation step, but we will not look into that here.
For Inverse \methodname, we use $S^{-1}B_2$ as initial guess as suggested by Proposition~\ref{t2}.
For the same reason, we use $A^{-1}B$ as initial guess for Projected \methodname.
Results for the MOC problem with $n_y=32$ are shown in \tabref{comp1}.
\llabel{nvec}Note that $B$ always has $n_y$ columns.
}

\begin{table}[ht]
\centering
\pink{
\small
\begin{tabular}{l|llllllll}
                        Method &  Rank &   Dim &   Its &   MVPs &  IMVPs & $t_s$ (s) &     $\rho$\\\hline
                   \methodname &    59 &   204 &   217 &   634 &     0 &     8 & $9.2\cdot 10^{-3}$\\
           Inverse \methodname &    60 &   127 &    34 &   128 &   128 &     7 & $9.5\cdot 10^{-3}$\\
         Projected \methodname &    80 &   133 &    36 &   134 &   134 &     9 & $9.9\cdot 10^{-3}$\\
                          EKSM &    60 &   448 &     7 &   224 &   256 &    11 & $6.2\cdot 10^{-3}$\\
                         PEKSM &    81 &   704 &    11 &   704 &   384 &    24 & $6.3\cdot 10^{-3}$\\
                          RKSM &    60 &   448 &    13 &   448 &   416 &    64 & $9.2\cdot 10^{-3}$\\
                         TRKSM &    60 &   187 &    16 &   219 &   155 &    46 & $9.7\cdot 10^{-3}$\\
            LR-ADI (heuristic) &    60 &   800 &    25 &     0 &   480 &   129 & $6.1\cdot 10^{-3}$\\
           LR-ADI (projection) &    60 &  1312 &    41 &     0 &   800 &   212 & $6.4\cdot 10^{-3}$
\end{tabular}
\caption{Comparison of different Lyapunov solvers.
  Rank is the rank of the final approximate solution, Dim is the maximum dimension of the approximation space during the iteration, Its is the amount of iterations, MVPs are the amount of matrix-vector products, IMVPs are the amount of inverse matrix-vector products and $t_s$ is the time required for solution of the Lyapunov equation, \blueii{which includes the computation of the LU-factorization when necessary.}
  For all methods the stopping criterion is a relative residual of $10^{-2}$.
  \methodname was restarted after 50 iterations with a tolerance of $4\cdot 10^{-6}$ for the eigenpairs that were retained.
  This is the same tolerance that was used for the other methods to minimize the rank of the approximate solution.
\label{tab:comp1}}
}
\end{table}\llabel{tab4}
\green{
Let us discuss all columns of \tabref{comp1} separately.
The first column contains the rank of the approximate solution, which is more-or-less the same for every method, because we did some post-processing on the non-\methodname methods to only keep eigenvectors belonging to eigenvalues that are larger than a certain tolerance.
In this case we took $4\cdot 10^{-6}$.
For \methodname, we do not have to do this since the restart method already takes care of this.
The projected variants have a larger rank, because they iterate on the full space instead of only the space belonging to the Schur complement.
}

\green{
Dim shows the maximum dimension of the search space during the iteration.
Note that this is much smaller for all variants of \methodname than for almost \pinkii{all of} the other methods.
For the standard \methodname method this is due to the restart that we use, but we restart after 50 iterations, so both other variants of \methodname do not restart, except in the last step when the method already converged and the rank is being minimized.
The only other method that has a small maximum space dimension is the Tangential Rational Krylov method, which also expands the space by only a few vectors in each iteration.
This is also the reason why we compare to this method.
}

\green{
\pinkii{The next column contains} the amount of iterations.
We show this just for reference purposes.
Every method has a different notion of what an iteration is, so we can not really compare these values.
For instance in this case \methodname expands by 3 vectors per iteration, Extended Krylov by 64, and LR-ADI by 32.
}

\green{
Then we show the MVPs, which are the amount of matrix-vector products.
For \methodname this is quite high compared to the other methods, but then the advantage is that no linear solves are \pinkii{required}, which can be seen if we look at the IMVPs.
If we include these in \methodname, which we did in the Inverse \methodname method, we see that it needs far fewer linear solves than the other methods.
}

\green{
The most important part of this table is the solution time $t_s$, from which we can see that all variants of \methodname are faster than all other methods.
The fastest method is Inverse \methodname.
However, it is only a bit faster than standard \methodname, at the cost of having to solve linear systems in each iteration.
For the solution of the linear systems, an LU-factorization was computed beforehand that was utilized by the two variants of \methodname that require a solve, and by the (Projected) Extended Krylov method.
\blueii{The time this takes is included in the solution time.}
For the LR-ADI and Rational Krylov methods precomputing the LU-factorization is not possible, because in each iteration a different shifted linear system has to be solved.
Mainly for this reason \methodname is much faster than Tangential Rational Krylov, even though it uses a larger space.
}

\green{
Lastly, we added a column with the explicitly computed final residual norms, which are all between $6\cdot 10^{-3}$ and $10^{-2}$.
}

\begin{figure}[ht]
  \centering\includegraphics[width=0.9\textwidth]{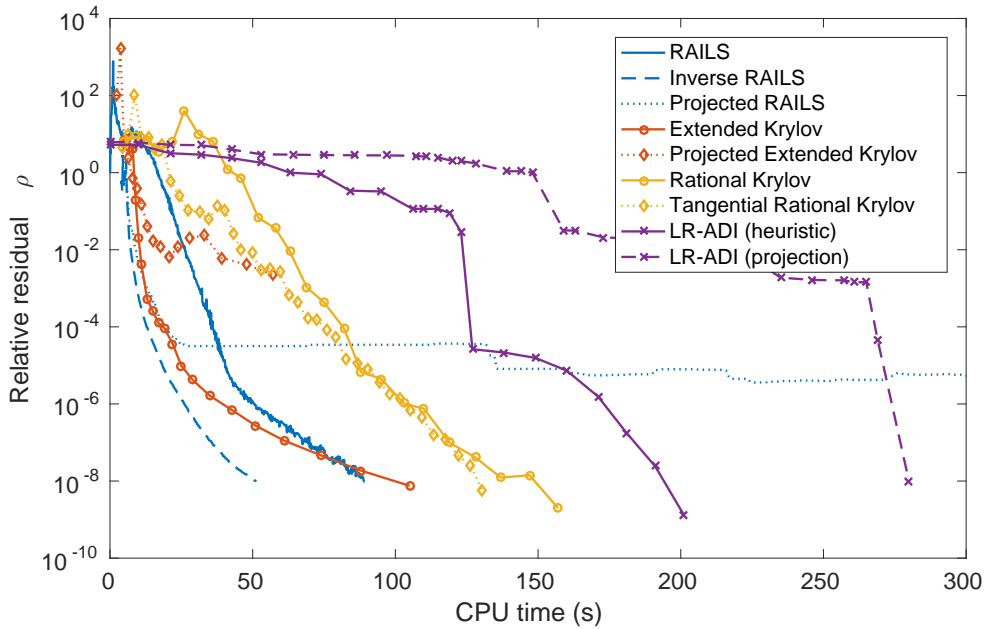}
  \caption{ \red {
Convergence history of the relative residual $\rho$ of all the different methods against CPU time.
\methodname was restarted after 15 iterations with a tolerance of $10^{-12}$ for the eigenpairs that were retained and expanded with 10 vectors per iteration.
}
 \label{fig:conv}}
\end{figure}

\red {
We use quite a loose tolerance for the results in \tabref{comp1}, which is sufficient for our purposes.
We will continue with this loose tolerance in later experiments.
However, with this tolerance, it is quite hard to see the convergence properties of the different methods.
\llabel{conv}In \figref{conv}, we show a convergence plot with a much smaller tolerance: a relative residual of $10^{-8}$.
We choose to put CPU time on the $x$-axis, since there is not really any other quantity that represents a similar amount of work for each method.
}

\begin{figure}[ht]
  \centering
\definecolor{mycolor1}{rgb}{0.49400,0.18400,0.55600}%
\definecolor{mycolor2}{rgb}{0.92900,0.69400,0.12500}%
\definecolor{mycolor3}{rgb}{0.85000,0.32500,0.09800}%
\definecolor{mycolor4}{rgb}{0.00000,0.44700,0.74100}%
\begin{tikzpicture}[font=\sffamily,scale=0.9]
\begin{axis}[%
width=4in,
height=2.82in,
scale only axis,
xmin=0,
xmax=2000,
xtick={0,500,1000,1500,2000},
xticklabels={{0},{500},{1000},{1500},{2000}},
xlabel={Maximum space dimension},
ymin=0,
ymax=8,
ytick={1,2,3,4,5,6,7},
yticklabels={{LR-ADI (projection)},{LR-ADI (heuristic)},{Tangential Rational Krylov},{Rational Krylov},{Extended Krylov},{Inverse RAILS},{RAILS}},
axis background/.style={fill=white}
]
\addplot[xbar,bar width=0.8,draw=black,fill=mycolor1,area legend] plot table[row sep=crcr] {%
1760	1\\
};
\addplot[xbar,bar width=0.8,draw=black,fill=mycolor1,area legend] plot table[row sep=crcr] {%
1248	2\\
};
\addplot[xbar,bar width=0.8,draw=black,fill=mycolor2,area legend] plot table[row sep=crcr] {%
448	3\\
};
\addplot[xbar,bar width=0.8,draw=black,fill=mycolor2,area legend] plot table[row sep=crcr] {%
832	4\\
};
\addplot[xbar,bar width=0.8,draw=black,fill=mycolor3,area legend] plot table[row sep=crcr] {%
1344	5\\
};
\addplot[xbar,bar width=0.8,draw=black,fill=mycolor4,area legend] plot table[row sep=crcr] {%
360	6\\
};
\addplot[xbar,bar width=0.8,draw=black,fill=mycolor4,area legend] plot table[row sep=crcr] {%
355	7\\
};
\end{axis}
\end{tikzpicture}%
  \caption{ \green {
Memory usage of the different methods as described in \figref{conv} in terms of maximum space dimension.
Both projected methods were left out of this figure since they did not converge.
}
 \label{fig:hbar}}
\end{figure}
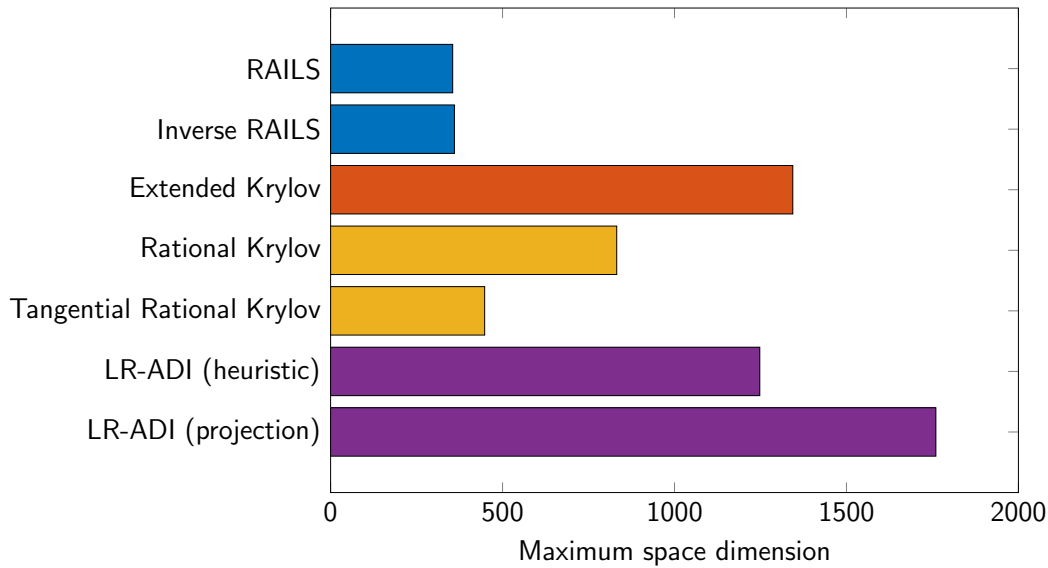

\red {
What we see in \figref{conv} is that two methods perform really badly: the Projected Extended Krylov method and Projected \methodname.
Projected Extended Krylov actually breaks down due to the small projected Lyapunov equation having eigenvalues with opposite signs and Projected \methodname stagnates.
This might be due to round-off errors during the projection, so it seems that here the projected variants are not very well suited for solving our problem.
The other methods all converge, but for the Krylov type methods, we clearly see that the cost per iteration increases the further they converge.
This is because it becomes increasingly expensive to solve the small projected Lyapunov equation.
On the other hand, the LR-ADI methods performed increasingly well for smaller tolerances.
They are still really expensive however, due to the shifted solves in every iteration.
The fastest method, again, is the Inverse \methodname.
It also shows almost monotone convergence, except in the first few steps, and converges very rapidly.
Standard \methodname performs very well, but convergence is more erratic.
It is, however, promising that a method that only uses matrix-vector products can converge faster than the other existing methods that we tried.
}

\green {
In \figref{hbar} we show the memory usage of the methods in \figref{conv}.
We left out the methods that did not converge.
Here it is clear that RAILS uses the least amount of memory due to the small amount of vectors that \pinkii{is} used to expand the space and due to the restart strategy.
}

\subsection{Numerical scalability} \label{sec:scal} \llabel{scal}

Now we want to show how \methodname behaves when we solve larger systems.
We do this by solving the same problem with different grid sizes, namely $4\times 32\times 16$, $4\times 64\times 16$ and $4\times 128\times 16$.
For the solution of the generalized Lyapunov equation, when increasing the dimension of our model by a factor 2, the size of the solution increases by a factor 4, since the solution is a square matrix.
\red {
However, we try to compute a low-rank approximation of the solution, so if the rank of the approximate solution stays the same when we increase the dimension of our model by a factor 2, the size required to store our approximate solution is also increased by a factor 2.
}

\red {
For the MOC problem, if we increase $n_y$ by a factor 2, we know that the amount of columns in $B$ also increases by a factor 2, since those two are equal.
From this, we would also expect the rank of the approximate solution to increase.
Therefore, we first look at a simplified problem where we use $\hat B = B \cdot \textbf{1}$ as right-hand side, so $\hat B$ is a single vector containing the row sums of the original $B$.
For our test problem, this can be seen as a fully space correlated stochastic forcing.
We expect that the rank of the approximate solution stays the same.
}

\red{
\llabel{t5}From \tabref{gridindep} we see that indeed the rank of the approximate solution does not change when we increase the dimension of the model.
However, the number of iterations to find the approximate solution is dependent on the spectral properties of $A$.
And since we are refining, new high frequency parts in the approximate solution will be amplified strongly  in the residual evaluation and appear also in the search space.
This will slow down the convergence process as can also be seen in \tabref{gridindep}.
}

\begin{table}[ht]
\centering
\green{
\begin{tabular}{l|llllll}
Size & Rank & Dim & Its & MVPs & $t_s$ (s) & $t_m$ (S)\\\hline
32   & 12   & 41  & 231 & 223  & 3  & 1.5\\
64   & 13   & 42  & 350 & 338  & 13 & 10 \\
128  & 13   & 42  & 536 & 518  & 58 & 52
\end{tabular}
\caption{Performance of \methodname for different grid sizes with $\hat B = B \cdot \textbf{1}$.
  The grids are of size $4 \times n_y \times 16$ where $n_y$ is the size in the first column.
  Rank is the rank of the final approximate solution, Dim is the maximum dimension of the approximation space during the iteration, Its is the amount of iterations, MVPs are the amount of matrix-vector products, $t_s$ is the time required for solution of the Lyapunov equation and $t_m$ is the cost of the matrix-vector product.
  The stopping criterion is a relative residual of $10^{-2}$.
  \methodname was restarted after 30 iterations with a tolerance of $10^{-5}$ for the eigenpairs that were retained.
  We expanded the space with 1 vector in each iteration.
\label{tab:gridindep}}
}
\end{table}

\red {
We are of course also interested in the increase of the solution time.
There are five things that influence this: the length of the vectors that span our basis, the amount of vectors that span our basis, the amount of iterations, the cost of the matrix-vector product, and the cost of solving the projected Lyapunov equation.
First of all, the length of the vectors that span our basis increases by a factor 2 when the dimension of the problem is increased by a factor 2, so we can expect a factor 2 in time increase from this.
Secondly, the rank of the approximate solution and the dimension of the search space does not increase.
This means that also the cost of solving the projected Lyapunov equation does not increase\pinkii{, so} we do not expect a solution time increase from this.
Then we have the amount of iterations, which does seem to increase by a factor 1.5.
And finally we have the cost of the matrix-vector product, which we listed in an extra column of \tabref{gridindep}.
Note that this is actually a matrix-vector product with the Schur complement $S$ which requires a linear solve with the matrix block $A_{11}$, which is relatively expensive.
If we subtract the cost of the matrix-vector product, we would expect a factor 3 in time cost increase from the increased vector length and the increased amount of iterations.
In \tabref{gridindep} we observe roughly a factor 2.
This being less than 3 might be due to higher efficiency of vector operations for larger vectors. 
}

\red {
Going back to the original problem, where we take $B$ as the original stochastic forcing, we expect the solution time to increase by a larger factor, since the projected Lyapunov equation solve and the vector operations become more costly, because the maximum search space dimension and the rank of the approximate solution now do increase.
In \tabref{gridsize} we see that the amount of iterations from $n_y=32$ to $n_y=64$ increases by a factor 1.5 again, so we would expect the solution time without matrix-vector products to increase by a factor 3, which is true.
The amount of iterations from $n_y=64$ to $n_y=128$ increases by a factor 2, so we would expect the time cost to increase by a factor 4, and this is also true.
This is because the solution of the projected Lyapunov equation is still relatively cheap for the problems we have here, since those projected equations have at most size $247\times 247$.
For problems with a larger maximum search space dimension, the cost of solving this projected Lyapunov equation would become dominant.
}
\begin{table}[ht]
\centering
\green{
\begin{tabular}{l|llllll}
Size & Rank & Dim & Its & MVPs & $t_s$ (s) & $t_m$ (s)\\\hline
32   & 59   & 198 & 233 & 682    & 7   & 2\\
64   & 86   & 223 & 340 & 997    & 31  & 17\\
128  & 147  & 247 & 652 & 1915   & 178 & 115
\end{tabular}
\caption{Performance of \methodname for different grid sizes.
  The grids are of size $4 \times n_y \times 16$ where $n_y$ is the size in the first column.
  Rank is the rank of the final approximate solution, Dim is the maximum dimension of the approximation space during the iteration, Its is the amount of iterations, MVPs are the amount of matrix-vector products and $t_s$ is the time required for solution of the Lyapunov equation and $t_m$ is the cost of the matrix-vector product.
  The stopping criterion is a relative residual of $10^{-2}$.
  \methodname was restarted after 50 iterations with a tolerance of $10^{-5}$ for the eigenpairs that were retained.
\label{tab:gridsize}}
}
\end{table}
\subsection{Towards a 3D model} \label{sec:3d} \llabel{3d}
\green{
Ultimately, we want to do computations on a full 3D model.
To illustrate what will happen for a full 3D model, we include some results where we take the forcing in such a way that it is not zonally averaged, but it is taken such that there is no correlation in the zonal direction.
}
\red {
\llabel{3db}The new forcing can be seen as a diagonal matrix, where all salinity nodes at the surface have a nonzero on the diagonal, and the rest of the matrix is zero.
We can write our new forcing, using Matlab notation, as $\hat B = \textrm{diag}(B \cdot \textbf{1})$ where $\hat B$ is the new forcing and $B$ is the original forcing.
This means that there are $4 \cdot n_y - 1$ nonzero columns in $\hat B$.
}

\green{
We choose to compare \methodname to Extended Krylov, which was the only method that came close to \methodname in the earlier experiments in terms of time, and to the Tangential Rational Krylov method, since this was the only method that came close in terms of maximum space dimension.
Since Extended Krylov does not perform well anymore when the projected system becomes really large, which would be the case for this problem, we split $\hat B$ into multiple parts.
We can do this because for our problem, the right-hand side $\hat B \hat B^{\top}$ can be written as
\begin{align*}
\hat B \hat B^{\top}=\sum_{i=1}^{4 \cdot n_y - 1} \hat B_i \hat B_i^{\top}
\end{align*}
where $\hat B_i$ is the $i$th column of $\hat B$.
Since the rank of the approximate solution is highly dependent on the amount of vectors in $\hat B$, we expect the maximum space dimension that is needed, and therefore also the size of the projected system, to decrease.
After splitting $\hat B$, we separately solve the Lyapunov equations with these new right-hand sides, of which we reduce the rank of the approximate solution separately to save memory.
Afterwards, we merge the low-rank solutions back into one low-rank solution, which is the approximate solution for the system with $\hat B$.
We report results with the amount of parts that resulted in the least amount of solution time.
The optimal amount of parts that we found was 4, so that is three of size $n_y$ and one of size $n_y-1$.
The maximum space dimension that we report is the maximum space size of solving one part plus the amount of vectors that are required to store the reduced approximate solution of the previous solves.
The results with the 3D forcing are shown in \tabref{grid3d}.
}

\begin{table}[ht]
\centering
\green{
\begin{tabular}{l|llllllll}
  Method                       & Size  & Rank  & Dim   & Its   & MVPs   & IMVPs  & $t_f$ (s) & $t_s$ (s)\\\hline
                          EKSM &   32  &   172 &   763 &     6 &  1114 &  1241 &        3  & 43\\
                               &   64  &   309 &  1353 &     5 &  1979 &  2234 &       15  &   278\\\hline
                         TRKSM &   32  &   177 &   681 &    16 &   808 &   554 & 0         &   109\\
                               &   64  &   314 &  1232 &    17 &  1487 &   977 & 0         &   645\\\hline
  \methodname                  &  32   &  166  &  312  &   252 &  739  & 0     & 0         & 16\\
                               &  64   &  276  &  419  &   406 &  1189 & 0     & 0         & 82\\
                               &  128  &  563  &  699  &   651 &  1912 & 0     & 0         & 584
\end{tabular}
\caption{Comparison of different Lyapunov solvers for different grid sizes with 
  $\hat B = \textrm{diag}(B \cdot \textbf{1})$.
  The grids are of size $4 \times n_y \times 16$ where $n_y$ is the size in the second column.
  Rank is the rank of the final approximate solution, Dim is the maximum dimension of the approximation space during the iteration, Its is the amount of iterations, MVPs are the amount of matrix-vector products, IMVPs are the amount of inverse matrix-vector products $t_f$ is the time required for computing the LU-factorization and $t_s$ is the time required for solution of the Lyapunov equation.
  For all methods the stopping criterion is a relative residual of $10^{-2}$.
  \methodname was restarted after 50 iterations with a tolerance of $10^{-6}$ for the eigenpairs that were retained.
  The tolerance that was used in Extended Krylov and Tangential Rational Krylov to minimize the rank of the approximate solution in such a way that the residual was still below the tolerance was set to $10^{-6}$ and $4\cdot 10^{-7}$ for the 32 and 64 problems respectively.
\label{tab:grid3d}}
}
\end{table}

\green {
We see that in this case \methodname can solve systems much faster and with much less memory than the Extended Krylov and Tangential Rational Krylov methods, even when applying the additional trick that we described above to reduce the memory usage of the Extended Krylov method.
It is also clear that both computing and applying the inverse becomes a real problem for these kinds of systems.
\methodname, however, does not need an inverse, and employs a restart strategy to reduce the space size, which is why it performs so much better.
The reason why we do not show any results with the Extended Krylov and Tangential Rational Krylov methods for the $4 \times 128 \times 16$ problem is that we did not have enough memory to do these computations, since we chose to not employ an iterative solver.
}
\section{Summary and Discussion}\label{sec:conc}

We have presented a new method \pinkii{for numerically determining covariance matrices, and hence we can compute probability density functions (PDFs) near steady states of deterministic PDE systems which are perturbed by noise.}
This method enables the application of the approach suggested by~\cite{Kuehn2011,Kuehn2015} to larger systems of SPDEs with algebraic constraints and \pinkii{exploits} the structure of the generalized Lyapunov equation in the computations.
This approach fits nicely within purely deterministic numerical bifurcation computations \cite[]{Keller1977} and methods to tackle the SPDEs directly \cite[]{Sapsis2009}.
It provides a Gaussian PDF which is valid under linearized dynamics near the deterministic steady state.

\green{
Our new algorithm to solve generalized Lyapunov equations is based on a projection method, where solutions are found iteratively by using a set of subspaces $V_k$.
The key new aspects are $(i)$ that at each iteration $k$, the subspace $V_k$ is expanded with the eigenvectors corresponding to the largest eigenvalues of the residual matrix $R$, orthogonalized with respect to $V_k$ and itself, and $(ii)$ the restart strategy that is used to keep the space dimension low.
We applied this method to a test problem consisting of a quasi two-dimensional model of the Atlantic Ocean circulation.
Our method, \methodname, outperforms both Krylov and ADI based methods \cite[]{Simoncini2007, Druskin2011, Stykel2012, Druskin2014, Kleinman1968, Penzl1999} for this test problem.
}

In practice, one has to provide the Jacobian matrix $A$ of a deterministic fixed point and the matrix $B$ representing the stochastic forcing in order to solve for the stationary covariance matrix $C$.
In a matrix-based pseudo-arclength continuation method the Jacobian is available since a Newton-Raphson method is used \cite[]{Dijkstra2014}.
The mass matrix $M$ is readily available from the model equations.
The method hence enables us to determine the continuation of local PDFs in parameter space.
In particular, the projection approach provides subspaces, which may be re-used directly or by computing predictors along continuation branches.

The availability of the new method opens several new directions of \pinkii{analysis}.
In one set of applications, $A$ is varied whereas the structure of $B$ is fixed.
This typically corresponds to varying a bifurcation parameter and analyzing the corresponding changes in PDF (i.e.~covariance matrix $C$) and transition probabilities (i.e.~overlap of PDFs of two steady states) as the steady states (i.e.~$A$) change and a bifurcation point is approached \cite[]{Kuehn2011}.
For example, this could be used in the test problem considered in this study in order to investigate the scaling law in PDF near the saddle-node bifurcation on the branch of pole-to-pole solutions.
In this way, the critical slowdown near a tipping point can be studied \cite[]{Mheen2013}.

For these types of problems, we expect \methodname to perform especially well since it can be restarted using the approximate solution of the previous continuation step, which should result in rapid convergence.

Another interesting application is to fix $A$ (or a set of $A$s that exists for fixed parameter values) and vary $B$.
In this case one steady state is fixed (or two in a regime of two steady states) and the impact of different $B$ (``noise products'') on the PDF/covariance matrix $C$ (and possibly transition probabilities) is investigated.
Different $B$ can be constructed by changing $(a)$ the dynamical component which is stochastically active (freshwater flux, wind, heat flux), $(b)$ the magnitude, and $(c)$ the pattern (i.e.~the cross-correlation structure and the spatial weighting of the auto-covariances).
Results from these computations will show the effect of the representation of small scale processes (the `noise') on the probability density function (under the restriction of linear dynamics).

For the ocean circulation problem it would be interesting to compare the different impacts of freshwater flux versus wind stress noise e.g.~on the PDF of the MOC or heat transports.
Moreover, regarding the wind stress and $(c)$ one could construct $B$ based on EOFs of the atmospheric variability, as for example considered in \cite{Dijkstra2008}.
Note that for studying the effect of wind-stress noise, a full three-dimensional model, as developed in \cite{DeNiet2007} is required.
In order to construct $B$ from more than one EOF one has to generalize the stochastic forcing to a sum of OU processes.
In that case we can use the linearity of the generalized Lyapunov equation and solve for each term separately and combine later on.

In future studies we aim to apply continuation methods, via Kramers'-type formulas adjoined to the continuation problem~\cite{Kuehn2011} as well as via transition path numerical methods~\cite{Henkelman2000}, to perform a full study of transition probabilities of transitions between MOC states in a three-dimensional global ocean model \cite[]{Thies2009}.
Knowledge of these transition probabilities is important to evaluate whether climate variability phenomena in the geological past (such as the Dansgaard-Oeschger events) could be caused by multiple steady states of the MOC and whether the transitions are stochastic or induced by crossing bifurcation points \cite[]{Ditlevsen2010}.
\vfill
\section*{Acknowledgements}
\green {
We would like to thank the creators of the M.E.S.S. package to make their code available.
We also thank Valeria Simoncini for making her code available and for giving suggestions on how to improve our results.
And finally we would like to thank the reviewers for their constructive comments.
}
This work is part of the Mathematics of Planet Earth research program, which is financed by the Netherlands Organization for Scientific Research (NWO).
The work of TEM and HAD was also carried out under the program of the Netherlands Earth System Science Centre (NESSC).
CK has been supported by an APART fellowship of the Austrian Academy of Sciences and by a Lichtenberg Professorship of the VolkswagenStiftung.

\clearpage 
\bibliographystyle{abbrvnatdoi}
\bibliography{JCP5_v1.0}

\end{document}